\documentclass[11pt]{article}
\usepackage{amsfonts,mathbbold,mathrsfs,
amsmath,amssymb,epsfig,subfigure,fancybox,graphicx,enumerate}
\newcommand{\wdh}{\widehat}

\newcommand{\ga}{\gamma}
\newcommand{\sg}{\sigma}

\def\l{\left|}
\def\r{\right|}

\newcommand{\wdt}{\widetilde}
\newcommand{\e}{\varepsilon}
\newcommand{\rr}{{\Bbb R}}
\newcommand{\M}{{\cal M}}
\newcommand{\cd}{(\cdot)}
\newcommand{\nd}{\noindent}
\newcommand{\tr}{{\rm tr}}

\def\para#1{\vskip .4\baselineskip\noindent{\bf #1}}
\def\qed{\strut\hfill $\Box$}

\newtheorem{thm}{Theorem}[section]
\newtheorem{prop}[thm]{Proposition}
\newtheorem{lem}[thm]{Lemma}
\newtheorem{cor}[thm]{Corollary}
\newtheorem{rem}[thm]{Remark}
\newtheorem{defn}[thm]{Definition}
\newtheorem{exm}[thm]{Example}

\newcommand{\thmref}[1]{Theorem~{\rm \ref{#1}}}
\newcommand{\lemref}[1]{Lemma~{\rm \ref{#1}}}
\newcommand{\corref}[1]{Corollary~{\rm \ref{#1}}}

\newcommand{\exmref}[1]{Example~{\rm \ref{#1}}}

\def\al{\alpha}

\newcommand{\K}{\wdh{\sf C}}

\newcommand{\beq}[1]{\begin{equation} \label{#1}}
\newcommand{\eeq}{\end{equation}}
\newcommand{\bed}{\begin{displaymath}}
\newcommand{\eed}{\end{displaymath}}
\newcommand{\bea}{\bed\begin{array}{rl}}
\newcommand{\eea}{\end{array}\eed}
\newcommand{\ad}{&\!\!\!\disp}
\newcommand{\aad}{&\disp}
\newcommand{\barray}{\begin{array}{ll}}
\newcommand{\earray}{\end{array}}
\def\({\left(}
\def\){\right)}
\def\disp{\displaystyle}

\numberwithin{equation}{section}

\begin{document}
\title{Stability of Nonlinear Regime-switching Jump Diffusions}
\author{Zhixin Yang,\thanks{Department of Mathematics,
Wayne State University, Detroit, Michigan 48202.  Email:
ef7538@wayne.edu. The research of this author was supported in
part by the National Science Foundation under DMS-0907753.}
\and
G. Yin\thanks{Department of Mathematics,
Wayne State University, Detroit, Michigan 48202.  Email:
gyin@math.wayne.edu. The research of this author was supported in
part by the National Science Foundation under DMS-0907753, and in
part by the Air Force Office of Scientific Research under
FA9550-10-1-0210.}}
\maketitle

\begin{abstract}
Motivated by networked systems, stochastic control, optimization,
and a wide variety of applications, this work is devoted to
systems of switching jump diffusions. Treating such  nonlinear
systems, we focus on stability issues. First asymptotic stability in
the large is obtained. Then the study on exponential $p$-stability
is carried out. Connection between almost surely exponential
stability and exponential $p$-stability is  exploited. Also
presented are smooth-dependence on the initial data. Using the
smooth-dependence, necessary conditions for exponential
$p$-stability are derived. Then criteria for asymptotic stability in
distribution are provided. A couple of examples are given to
illustrate our results.

\vskip 0.1 true in \nd{\bf Key Words.} jump diffusion, switching
process, stability in the large,
smooth dependence on initial data, stability in distribution.

\vskip 0.1 true in \nd{\bf Brief Title.} Stability of Jump Diffusions
with Switching
\end{abstract}

\newpage
\setlength{\baselineskip}{0.22in}
\section{Introduction}
Randomly varying switching  systems have drawn
increasing attention
recently, especially in
the fields of control and optimization.
This is largely
owing to their ability to model complex systems, which
can be used in a wide range of applications
in
consensus controls, distributed computing,
autonomous or semi-autonomous vehicles,
multi-agent systems, tele-medicine, smart grids, and
financial engineering etc.
One of the common features of the many systems mentioned above is
that they
may be represented as
networked systems.
In a typical networked system, different nodes are connected through
a communication link described by a network topology or
configuration. Most work to date dealt with networked systems with
fixed topology. Nevertheless, data routing, packet aggregations,
channel uncertainties, and switching links to different network hubs
demand the consideration of topology changes in a networked system.
Thus fixed topology becomes inadequate and random
environment and uncertainty must be taken into consideration.
For example, in the original formulation
of consensus problems \cite{Rey,JT,VCBCS},
one dealt with a fixed configuration or topology, whereas
consideration of randomly varying topologies leads to switching
diffusion processes
\cite{YinSW}.

Facing the demands and pressing needs, this paper considers
systems that are formulated as
regime-switching jump
diffusions.
Because many systems in networked systems are in operation for very
long time, their asymptotic behavior, namely, stability is of crucial
importance;
see \cite{BadowskiY,Liu}  and references therein for related work.
Due to the involvement of multiple
stochastic processes, care must be taken to treat the stability
issues, which is the objective of the current paper.

One of the main features of the
 underlying systems we consider here is the coexistence of
 dynamics described
 by solutions of differential equations
 and discrete events whose values belong to a finite
 set; see \cite{Hespanha05,KarM} and
 references therein.
The usual formulation in the traditional dynamic system setup
described by differential or difference equations alone becomes
unsuitable. A class of models naturally replacing the traditional
setup is a process with two components in which one of them
delineates the
 dynamics that may be represented
as a solution
of a differential equation and the other portraits the discrete event
movements (see \cite{BW87} for an example in finance application).
To take into consideration of possible inclusion of the
 Poisson type of random
processes, we consider jump diffusion processes with random
switching.

In recent years, switching stochastic systems
 have  received much attention; see
 \cite{Mao199945,MaoY06,Maoyuan} and references therein for a
systematic treatment on Markov modulated switching
diffusions; see also
\cite{Maoyuan1} for
stability of switching diffusions with delays. In addition, switching
diffusion with continuous dependence on initial data
were treated in
\cite{YinZhu}.
Concerning jump diffusions, we refer the reader to
\cite{LiDS,MR99,Wee} for the study on such properties
as ergodicity and
stability. Switching jump diffusions with state dependent
switching have also been examined in \cite{Xi,Yinxi,YinZ10} etc.,
in which stability in probability, asymptotic
stability in probability, and almost surely exponential stability
were dealt with.
Our aims in this paper are to establish
a number of
results on different modes of stability
that have not been studied for switching jump diffusions
to date to the best of our knowledge. We begin with
asymptotic stability in the large,
proceed to exponential
$p$-stability and obtain smooth dependence on the initial data. As a
nice application of the smooth dependence on the initial data, we
derive necessary conditions for $p$-stability, which can be viewed
as a Lyapunov converse theorem. The aforementioned results all begin
with an equilibrium point of the switching jump diffusion. In
absence of information of the equilibrium, an appropriate notion of
stability is stability in distribution. Under simple conditions, we
obtain sufficient conditions for asymptotic stability in
distribution.

The rest of the paper is arranged as follows. We begin with the
precise description of the system in Section \ref{sec:form}. Section
\ref{sec:as-lar} concentrates on asymptotic stability in the
large. Section \ref{sec:ex-p-s} proceeds to the study on exponential
$p$-stability. Section \ref{sec:Lya} furthers our investigation by
examining the smooth-dependence on the initial data, and Section
\ref{sec:as-d} presents criteria for asymptotic stability in
distribution. Section \ref{sec:exm} presents a few examples to
demonstrate our results. Finally, the paper is concluded in section
\ref{sec:fur} with further remarks.

\section{Formulation}\label{sec:form}

This section presents the formulation of the problem. We begin with
certain notation needed together with a number of definitions.
 We use $z'$ to denote the transpose of $z\in \mathbb{R}^{l_1
\times l_2}$ with $l_i \ge 1$, and $\mathbb{R}^{r\times 1}$ is
simply written as $\mathbb{R}^r$. Denote
$\mathbbold{1}=(1,1,\dots,1)'\in
\mathbb{R}^r$ that is a column vector with all entries being 1. For a
matrix $A$, its trace norm is denoted by $|A|=\sqrt{\tr(AA')}$. Let
$\left(X(t),\alpha(t)\right)$ be a two-component Markov process in
which $X(\cdot)$ takes values in $\mathbb{R}^r$
 and $\al\cd$ is a switching process taking values in a finite set
$\M=\{1,2,3,\dots,m\}$. Let $\Gamma$ be a subset of $\rr^r -\{0\}$
that is the range space of the impulsive jumps. For any subset $B$
in $\Gamma$, $N(t,B)$ counts the number of impulses on $[0,t]$ with
values in $B$; $b(\cdot,\cdot): \rr^r \times \M \mapsto \rr^r$,
$\sg(\cdot,\cdot): \rr^r \times  \M
 \mapsto \rr^r \times \rr^d$, and
$g(\cdot,\cdot,\cdot): \rr^r \times \M \times \Gamma \mapsto \rr^r$
are suitable functions whose precise conditions will be given later.
Consider the dynamic system given by \beq{2.1} \barray \ad
dX(t)=b(X(t),\alpha(t))dt+\sigma(X(t),\alpha(t))dw(t)+dJ(t),
\\
\ad  J(t)=\int^t_0\int_\Gamma
g(X(s^-),\alpha(s^-),\gamma)N(ds,d\gamma),\\
\ad X(0)=x,\,\alpha(0)=\alpha,\earray \eeq where
the switching process $\al(\cdot)$
obeys
the transition rule
\beq{2.11}
P\{\alpha(t+\Delta t)=j|\alpha(t)=i,
X(s),\alpha(s), s\leq t
\}=q_{ij}(X(t))\Delta t+o(\Delta t),\text{ for } i\neq j,
\eeq
 $w(t)$ is a
$d$-dimensional standard Brownian motion, and $N(\cdot,\cdot)$ is a
Poisson measure such that the jump process $N(\cdot,\cdot)$ is independent of the
Brownian motion $w(\cdot)$.
Equation \eqref{2.1} can be written as integral form:
\[X(t)=x+\int^t_0b(X(s),\alpha(s))ds+\int^t_0\sigma(X(s),
\alpha(s))dw(s)+\int^t_0\int_\Gamma g(X(s^-),
\alpha(s^-),\gamma)N(ds,d\gamma).\] Here we have used a setup
similar to \cite{Yinxi}. When we wish
to emphasize the initial data
 dependence in the sequel, we write the
process as $(X^{x,\alpha}(t),\alpha^{x,\alpha}(t))$.
Note that
although the two-component process $(X(t),\al(t))$ is Markov,
$\al(t)$ generally is not a Markov chain due to the dependence
of the state $x$ in the generator.
The transition rule indicates that $\al(t)$ depends on the
jump diffusion component.
 Thus the setup we consider is more
general than that of considered in the literature, whereas in the past work
it was often assumed that $\al(t)$ itself is a Markov chain and
$w(t)$ and $\al(t)$ are independent.

For future use, we define a compensated or centered Poisson measure
as
\[\wdt {N}(t,B)=N(t,B)-\lambda t\pi(B)\text{  for  }B\subset\Gamma,\]
where $0<\lambda<\infty$ is known as the jump rate and $\pi(\cdot)$
is the jump distribution (a probability measure). With this centered
Poisson measure, we can rewrite $J(t)$ as
\[J(t)=\int^t_0\int_\Gamma g(X(s^-),
\alpha(s^-),\gamma)\wdt {N}(ds,d\gamma)+ \lambda\int^t_0\int_\Gamma
g(X(s^-),\alpha(s^-),\gamma)\pi(d\gamma)ds,\] which is the sum of a
martingale and an absolute continuous process provided certain
conditions are satisfied for the function $g(\cdot)$.

Note that the evolution of the discrete component $\alpha(\cdot)$
can be represented by a stochastic integral with respect to a
Poisson measure (e.g., \cite{Skoro}). For $x\in\rr^{r}$ and $i,j\in
\mathcal{M}$ with $j\neq i$, let $\Delta_{ij}(x)$ be the consecutive
(with respect to the lexicographic ordering on $\mathcal{M}\times
\mathcal{M}$), left-closed, right-open intervals of the real line,
each having length $q_{ij}(x)$. Define a function $h:\rr^r \times \M
\times \mathbb{R} \mapsto \mathbb{R}$ by \beq{2.3}
h(x,i,z)=\sum^{m}_{j=1}(j-i)I_{\{z\in \Delta_{ij}(x)\}}. \eeq That
is, with the partition $\{\Delta_{ij}(x):i,j\in \mathcal{M}\}$ used
and for each $i\in\mathcal{M}$, if $z\in \Delta_{ij}(x)$,
$h(x,i,z)=j-i$; otherwise $h(x,i,z)=0$. Then we may write the
switching process as a stochastic integral
 \beq{2.4}
d\alpha(t)=\int_{\mathbb{R}}h(X(t),\alpha(t-),z)N_1(dt,dz), \eeq
where $N_1(dt,dz)$ is a Poisson random measure with intensity
$dt\times \wdt m(dz)$, and $\wdt m(\cdot)$
is the Lebesgue measure on
$\mathbb{R}$. The Poisson random measure $N_1(\cdot,\cdot)$ is
independent of the Brownian motion $w(\cdot)$ and the Poisson
measure $N(\cdot,\cdot)$. For subsequent use, we define another
centered Poisson measure as
 \bea
\mu(dt,dz)=N_1(dt,dz)-dt\times\wdt m(dz). \eea

The generator $\mathcal {G}$ associated with the process
$(X(t),\al(t))$ is defined as follows: For each $i \in \mathcal
{M}$, and for any twice continuously differentiable function
$f(\cdot,i)$, \bea \ad \mathcal {G} f(x,\cdot)(i)=\mathcal
{L}f(x,\cdot)(i)+\lambda \int_\Gamma[f(x+g
(x,i,\gamma),i)-f(x,i)]\pi (d\gamma),\eea
 where $\mathcal{L}$ is the
operator for a switching diffusion process  given by
\beq{operator-def} \barray
 {\cal L}f(x,\cdot)(i)\ad= \frac{1}{2}\sum_{k,l = 1}^r {a_{kl}
(x,i)\frac{{\partial ^2 f(x,i)}}{{\partial x_k \partial x_l }}}  +
\sum_{k = 1}^r {b_k (x,i)\frac{{\partial f(x,i)}}{{\partial x_k }}}  +
Q(x)f(x,\cdot)(i) \\
 \ad = \frac{1}{2}\tr(a(x,i)Hf(x,i))
 + b'(x,i)\nabla f(x,i) + Q(x)f(x,\cdot)(i),
 \quad i \in {\cal M},
 \earray\eeq
where $x\in\rr^r$, $a(x,i)=\sg(x,i)\sg'(x,i)$, $\nabla f(\cdot, i)$
and $Hf(\cdot, i)$ denote the gradient and Hessian matrix of
$f(\cdot, i)$, respectively, and
$Q(x)=(q_{ij}(x))$ is an $m\times m$
matrix depending on $x$ satisfying the
 $q$-property, namely,
(i) $q_{ij}(x)$ is Borel measurable and uniformly bounded for all
$i,j\in \mathcal{M}$ and $x\in \mathbb{R}^r$;

(ii) $q_{ij}(x)\ge 0$ for all $x\in \rr^r$ and $i\neq j;$

(iii) $q_{ii}(x)=-\sum \limits_{j\ne i\atop j \in {\cal M}}
q_{ij}(x)$ for all $x\in \rr^r$ and $i\in \mathcal{M}$. Denote

\bea Q(x)f(x,\cdot)(i)=\sum\limits_{j\in
\M}q_{ij}(x)f(x,j)=\sum\limits_{j\neq i\atop j\in
\cal{M}}q_{ij}(x)(f(x,j)-f(x,i)),i\in \mathcal{M}.\eea In what
follows, we often write $\mathcal{L}f(x,\cdot)(i)$ as
$\mathcal{L}f(x,i)$ and $\mathcal{G}f(x,\cdot)(i)$ as
$\mathcal{G}f(x,i)$ for convenience whenever there is no confusion.
By virtue of the generalized It\^{o}'s formula, we have that
\[f(X(t),\alpha(t))-f(x,\alpha)-\int^t_0 \mathcal
{G}f(X(s),\alpha(s))ds \ \hbox{ is a martingale.}\]

To proceed, we need the following assumptions.
\begin{enumerate}[({A}1)]
  \item  The functions $b(\cdot,i)$, $\sigma(\cdot,i)$,
  and $g(\cdot,i,\ga)$ satisfy
  $b(0,i)=0$, $\sigma(0,i)=0$,
  and $g(0,i,\ga)=0$ for each $i\in\M$; $\sigma(x,i)$ vanishes only
at $x=0$
  for each $i\in\M$.

  \item There exists a positive constant
  $K_0$ such that for each $i\in\mathcal {M},x,y \in
  \rr^r$ and $\gamma \in \Gamma$,
\bea \ad |b(x,i)-b(y,i)| +|\sigma(x,i)-\sigma(y,i)|\leq K_0|x-y|,\\
\ad |g(x,i,\gamma)-g(y,i,\gamma)|\leq K_0|x-y|.\eea

  \item There exists $g^*(i)$ satisfying
  \[|g(x,i,\gamma)|\le g^*(i)|x|\text{  for each  }
  x\in\mathbb{R}^r, i\in\M,
  \text{  and each  }\gamma\in\Gamma.\]
\end{enumerate}

We elaborate on the conditions briefly. Condition (A1) indicates that
$0$ is an equilibrium point; (A2) is a Lipschitz condition on the
functions. It together with the equilibrium point 0 implies that the
functions grow at most linearly.
Several of our results to follow are concerned with equilibrium
point of the switching jump diffusions. To proceed, we make the
following definitions by adopting the terminologies of
\cite{YinZ10}.

\begin{defn}\label{def-st}
{\rm The equilibrium point $x=0$ of system \eqref{2.1} and
\eqref{2.11} is said to be

\begin{enumerate}[(i)]

\item {\em stable in probability},
if for any $ \e  >0$ and any $\al\in\M$, $ \lim\limits_{x\to 0}
P\{\sup\limits_{t\ge 0} |X^{x,\al}(t)| > \e \} =0$; and $x=0$ is
said to be {\em unstable in probability} if it is not stable in
probability.

\item {\em asymptotically stable in probability}, if
it is stable in probability and $ \lim\limits_{x\to 0}
P\{\lim\limits_{t\to \infty}X^{x,\al}(t)=0\}=1, \hbox{ for any
}\al\in\M;$

  \item {\em asymptotically
stable in the large}, if it is stable
  in probability and
$P\{\lim\limits_{t\rightarrow \infty}
 X^{x,\alpha}(t)=0\}=1,
 \text{  for any  }(x,\alpha)\in \rr^r \times \M;$

  \item
{\em  exponential
  $p$-stable}, if
  for some positive constants $K$ and $k$,
$E|X^{x,\alpha}(t)|^p \leq K|x|^pe^{-kt}$,  for any $(x,\alpha)\in
\rr^r \times \M;$

\item {\em almost surely exponential
stable}, if for any $(x,\alpha)\in\mathbb{R}^r\times\cal{M}$, $
\limsup\limits_{t \to \infty } \frac{1}{t}\ln (|X^{x,\alpha}(t)|) <
0$ w.p.1.
\end{enumerate}
}\end{defn}

As a preparation, we first recall a lemma, which indicates that the
equilibrium $(0,\al)$ is inaccessible in that starting with any
$x\not =0$, the system will not reach the origin with probability
one. The proof of this lemma can be found in \cite[Lemma
2.10]{Yinxi}.

\begin{lem} \label{0-inacc}
$P\{X^{x,\al}(t)\not= 0, t\ge 0\}=1$, for any $x\not=0$
and $\al\in\M$.
\end{lem}

\section{Asymptotic Stability in the Large}
\label{sec:as-lar} To proceed, we first recall two lemmas. The
detailed proof can be found in \cite{Yinxi}.

\begin{lem}\label{lem3.1}
  Let $D\subset \mathbb{R}^r$ is a neighborhood of $0$. Suppose that
for each $i\in \cal{M}$, there exists a nonnegative Lyapunov
function $V(\cdot,i):\,D\mapsto \mathbb{R}$ such that
\begin{itemize}
  \item[{\rm(i)}] $V(\cdot,i)$ is continuous in $D$ and vanishes only at $x=0$;
  \item[{\rm (ii)}] $V(\cdot,i)$ is twice continuously differentiable in
  $D-\{0\}$ and satisfies
  $ \mathcal {G}V(x,i)\leq 0$   for all
$x\in \ D-\{0\}$.
\end{itemize}
Then the equilibrium point $x=0$ is stable in probability.
\end{lem}

Define \beq{3.1}
 \tau_{\rho, \varsigma}:=\inf\{t\geq 0:|X(t)|=\rho
\hbox{  or  } |X(t)|=\varsigma\}, \eeq for any $0< \rho < \varsigma$
and any $(x,\alpha)\in \mathbb{R}^r \times \M$ with $\rho < |x| <
\varsigma$.

\begin{lem}\label{lem3.2}
Assume that the conditions of \lemref{lem3.1} hold, and that for any
sufficiently small $0<\rho<\varsigma$ and any
$(x,\alpha)\in\mathbb{R}^r \times \M$ with $\rho < |x| < \varsigma$,
$ P\{\tau_{\rho,\varsigma}<\infty\}=1.$ Then the equilibrium point
$x=0$ is asymptotically stable in probability.
\end{lem}

\begin{thm}\label{thm:st-l}
Assume that the conditions of \lemref{lem3.2} hold, and that
$V_\varsigma: =\inf\limits_{|x|\ge \varsigma \atop i\in \mathcal
{M}} V(x,i)\rightarrow \infty$  as $\varsigma\rightarrow \infty$.
 Then the equilibrium point $x=0$ is asymptotically stable in the
large.
\end{thm}

\para{Proof.} For each $i \in \mathcal{M}$,
for any $\varepsilon > 0$ and $(x,\alpha) \in \mathbb{R}^r\times
\cal{M}$, there exists a
 $\varsigma >|x|$ large enough such that
 $\inf\limits_{|X|\ge \varsigma
 \atop i\in\M}V(X,i)\ge 2V(x,\alpha)/{\varepsilon}$.

 Let $\tau_{\varsigma}$ be
 the stopping time $\tau_{\varsigma}:=\inf\{t\ge
0:|X(t)|\ge
 \varsigma\}$ and $t_\varsigma=\tau_\varsigma\wedge t$. Then it
 follows from Dynkin's formula that
 \[EV(X(t_\varsigma),\alpha(t_\varsigma))-
V(x,\alpha)=E\int_0^{t_\varsigma}\mathcal {G}V(X(u),\alpha(u))du\le
0.\]
 Consequently,
 $EV(X(t_\varsigma),\alpha(t_\varsigma))\le V(x,\alpha).$
 Then we have
 \[E[V(X(\tau{_\varsigma)},
 \alpha(\tau{_\varsigma}))I_{\{\tau_{\varsigma} < t\}}]
 \le V(x,\alpha).\]
 Hence,
 $\frac{2V(x,\alpha)}{\varepsilon}
 P(\tau_{\varsigma} < t)\le V(x,\alpha).$
So
 $P(\tau_{\varsigma} < t)\le \varepsilon/2.$
 Let $t\rightarrow \infty $,
 $P(\tau_{\varsigma}<\infty)\le\varepsilon/2.$
Then it follows from \lemref{lem3.2} that, for any $\rho>0$ with
$\rho<|x|<\varsigma$ we have
\[1=P(\tau_{\rho,\varsigma}<\infty)\le
P(\tau_\rho<\infty)+P(\tau_{\varsigma}<\infty),\] in which
$\tau_\rho$ is the stopping time $\tau_\rho:=\inf\{t\ge 0:|X(t)|\le
\rho\}$, where $\tau_{\rho,\varsigma}$ was defined in \eqref{3.1}.
 Consequently,
$P(\tau_\rho <\infty)\ge 1-\varepsilon/2.$ This implies that
$P\{\mathop {\inf\limits_{t \ge 0} |X(t)|} \le \rho \} \ge 1 -
\varepsilon/2.$ Since $\rho>0$ can be arbitrarily small, $P\{\mathop
{\inf\limits_{t \ge 0}|X(t)|}   = 0\} \ge 1 - \varepsilon/2.$

Now we can follow the same techniques in \cite[Lemma 7.6]{YinZ10}
and obtain $P\{\lim\limits_{t\to \infty}X(t)=0\}\ge1-\varepsilon/2.$
That is, the equilibrium point $x=0$ is asymptotically stable in the
large as desired. \qed

For application, it is important to be able to handle linearized
systems.
Similar to \cite{Yinxi}, we pose the following condition.

\begin{itemize}
\item[(A4)]
 For each $i\in
\M$, there exist $b(i),\sigma_l(i)\in \rr^{r\times r}$ for
$l=1,2,\dots,d$, and a generator of a continuous-time Markov chain
$\widehat{Q}=(\hat{q}_{ij})$ with the corresponding Markov chain
denoted by $\wdh \al(t)$ such that as $x\rightarrow 0$,
 \bea
  b(x,i)\ad=b(i)x+o(|x|),\\
\sigma(x,i)\ad=(\sigma_1(i)x,\sigma_2(i)x,...,\sigma_d(i)x)+o(|x|),\\
Q(x)\ad=\widehat{Q}+o(1) .\eea
 Moreover, $\widehat{Q}$ is irreducible.

\end{itemize}

Assumption (A4) indicates that near the origin, the coefficients are
locally linear.
By choosing a Lyapunov function properly, we have the same
sufficient condition for asymptotically stable in the large as that of
asymptotically stable in probability. The result is provided below,
and the proof is omitted. The method involved is similar to
\cite[Theorem 3.5]{Yinxi}.

\begin{cor}\label{cor-34}
Under assumptions {\rm(A1)}-{\rm (A4)}, the equilibrium point $x=0$
of the system given by \eqref{2.1} and \eqref{2.11} is
asymptotically
stable in the large if
\[
\sum_{i \in \M} {\mu _i \left(\Lambda _{\max} (\frac{{b(i) +
b'(i)}}{2}) + \frac{1}{2}\Lambda _{\max} \left(\sum_{l = 1}^d
{\sigma _l (i)\sigma' _{l} (i)} \right) + \lambda g^* (i)\right) <
0}.
\]
\end{cor}

In which $\mu=(\mu_1, \mu_2, \cdots, \mu_m)\in \mathbb{R}^{1\times
m}$ is the stationary distribution of $\wdh \al(t)$ and
$\Lambda_{\max}(A)$ denotes the largest eigenvalue of the symmetric
part of $A$.

\section{Exponential $p$-stability}\label{sec:ex-p-s}
In this section, we give a sufficient condition for exponential
$p$-stability. To proceed, we first recall a lemma, which indicates
that the process $\left(X(t),\alpha(t)\right)$ has no finite
explosion time, also known as regular. The proof of this lemma can
be found in \cite[Lemma 2.8]{Yinxi}.

\begin{lem}\label{lem41}
  Under assumptions {\rm(A1)}-{\rm(A3)}, the switching jump
  diffusion $(X(t),\alpha(t))$ is regular.
\end{lem}

\begin{thm}\label{Theorem 4.1} Let $D\subset \mathbb{R}^r$ be a
neighborhood of $0$. Assume that the conditions of \lemref{lem41}
hold and assume that for each $i\in\mathcal{M}$, there exists a
nonnegative Lyapunov function $V(\cdot,i)$: $D\mapsto\mathbb{R}$
such that $V(\cdot,i)$ is twice continuously differentiable in
$D-\{0\}$, and satisfies the following conditions:
 \beq{4.1}k_1|x|^p\le V(x,i)\le k_2|x|^p,\quad
x\in D,\eeq
 \beq{4.2}\mathcal{G}V(x,i) \le-kV(x,i)\text{ for all
}x\in D-\{0\},\eeq
 for some positive constants $k_1, k_2$ and $k$. Then the equilibrium point
$x=0$ is exponential $p$-stable.
\end{thm}

\para{Proof.}
Consider a sequence of  stopping times $\{\tau_n\}$ defined by
$\tau_n:=\inf\{t\ge 0: |X(t)|\ge n\}$ and let $t_n=t\wedge \tau_n$.
By virtue of Dynkin's formula and \eqref{4.2}, since
$kV(X(s),\alpha(s)) +\mathcal{G}V(X(s),\alpha(s))\le 0,$ so
\bea\disp E[e^{k(t_n)}V(X(t_n),\alpha(t_n))]
\ad=V(x,\alpha)+E\int^{t_n}_0[e^{ks}(kV(X(s),\alpha(s))
+\mathcal{G}V(X(s),\alpha(s)))]ds\\
\ad\le V(x,\alpha).\eea Let $n\rightarrow\infty$, by Fatou's Lemma
and  \lemref{lem41}, we have
 $ E[e^{kt}V(X(t),\alpha(t))]\le V(x,\alpha).$
 Hence,\bea \ad e^{kt}E(k_1|X(t)|^p) \le e^{kt}EV(X(t),\alpha(t)) \le
V(x,\alpha)\le
 k_2|x|^p.\eea
Then we obtain $E|X(t)|^p\le K|x|^pe^{-kt}$. The theorem is thus
proved. \qed

\begin{rem}\label{Rem1}{\rm
 In the above and hereafter, $K$ is used as a generic
positive constant, whose value may be different in
 different appearances.
Under the
  conditions of \thmref{Theorem 4.1}, we can also obtain the result of
  almost surely exponential stability by
  similar argument in \cite[Theorem 5.8.1]{RZ}.}\end{rem}

\begin{thm}
\label{Theorem 4.2} Under assumptions {\rm(A1)}-{\rm(A3)},
exponential $p$-stability implies almost surely exponential
stability.
\end{thm}

\para{Proof.}
Because
\[X(t)=x+\int^t_0b(X(s),\alpha(s))ds
+\int^t_0\sigma(X(s),\alpha(s))dw(s)+\int^t_0\int_\Gamma
g(X(s^-),\alpha(s^-),\gamma)N(ds,d\gamma),\]
 we have for any $p\ge 2$ that
\bea |X(t)|^p\ad\le 4^{p-1}[|x|^p+\left|\int^t_0
b(X(s),\alpha(s))ds\right|^p+\left|\int^t_0
\sigma(X(s),\alpha(s))dw(s)\right|^p\\
\aad \qquad +\left|\int^t_0\int_\Gamma
g(X(s^-),\alpha(s^-),\gamma)N(ds,d\gamma)\right|^p].\eea
For any $t$, there exists such an $n$ that $t\in [n-1,n]$ and the
following inequality holds
 \beq{4.3}\barray
 \ad
E\left[\sup_{n-1\le t\le n}|X(t)|^p\right]\\
\aad \ \le 4^{p-1} E|X(n-1)|^p+4^{p-1}E\left(\sup_{n-1\le
t\le
n}\left|\int^{t}_{n-1}b(X(s),\alpha(s))ds\right|^p\right)\\
\aad\quad +4^{p-1}E\left(\sup_{n-1\le t\le
n}\left|\int^{t}_{n-1}\sigma(X(s),\alpha(s))dw(s)\right|^p\right)\\
\aad \quad +4^{p-1}E\left(\sup_{n-1\le t\le
n}\left|\int^{t}_{n-1}\int_\Gamma
g(X(s^-),\alpha(s^-),\gamma)N(ds,d\gamma)\right|^p \right).
\earray\eeq
By (A1) together with the H\"older inequality and the martingale
inequality,
\bea \ad E\left(\sup_{n-1\le t\le n}
\left|\int^{t}_{n-1}b(X(s),\alpha(s))ds\right|^p \right)\le
K\int^n_{n-1}E|X(s)|^pds,\\
\ad E\left[\sup_{n-1\le t\le
n}\left|\int^{t}_{n-1}\sigma(X(s),\alpha(s))dw(s)
\right|^p\right]\le K\int^n_{n-1}E|X(s)|^pds. \eea

For the Poisson jump part,
\beq{4.4}
 \barray
 \ad E\left[\sup_{n-1\le t\le
n}\left|\int^{t}_{n-1}\int_\Gamma
g(X(s^-),\alpha(s^-),\gamma)N(ds,d\gamma)\right|^p\right]\\
\aad \ =
 E\Big[\sup_{n-1\le t\le
n}\Big|\int^{t}_{n-1}\int_\Gamma g(X(s^-),\alpha(s^-),\gamma)
\wdt {N}(ds,d\gamma)\\
\aad \qquad +\lambda\int^{t}_{n-1}\int_\Gamma
g(X(s^-),\alpha(s^-),\gamma)\pi(d\gamma)ds\Big|^p\Big]\\
\aad\ \le
 2^{p-1}E\left[\sup_{n-1\le t\le
n}\left|\int^{t}_{n-1}\int_\Gamma
g(X(s^-),\alpha(s^-),\gamma)\wdt {N}(ds,d\gamma)\right|^p\right]\\
\aad\qquad+2^{p-1}E\left[\lambda\sup_{n-1\le t\le
n}\left|\int^{t}_{n-1}\int_\Gamma
g(X(s^-),\alpha(s^-),\gamma)\pi(d\gamma)ds\right|^p \right]. \earray
\eeq

Using H\"{o}lder inequality and assumptions {\rm(A1)}-{\rm(A3)} for
the last term of \eqref{4.4},
 detailed
computation leads to
 \bea \disp \ad
E\left[\lambda\sup_{n-1\le t\le
n}\left|\int^{t}_{n-1}\int_\Gamma
g(X(s^-),\alpha(s^-),\gamma)\pi(d\gamma)ds\right|^p\right]\\
\aad \ \le KE\left[\int^{n}_{n-1}\int_\Gamma
|X(s^-)|^p\pi(d\gamma)ds\right]\\
\aad\ \le K\int^n_{n-1}E|X(s^-)|^pds. \eea Now let us handle the
martingale part in \eqref{4.4}. By H\"{o}lder inequality,
assumptions {\rm(A1)}-{\rm(A3)}, and properties of stochastic
integral with respect to a Poisson measure, we have \bea \ad
E\left[\sup_{n-1\le t\le n}\left|\int^{t}_{n-1}\int_\Gamma
g(X(s^-),\alpha(s^-),\gamma)\wdt {N}(ds,d\gamma)\right|^p\right]
\\
\aad \ \le K E\left( \int^{n}_{n-1}\int_\Gamma\left|
g(X(s^-),\alpha(s^-),\gamma)\right|^2ds\pi (d\gamma)\right)^{p/2}\\
\aad \ \le KE\left( \int^{n}_{n-1}\int_\Gamma
|X(s^-)|^2ds\pi (d\gamma)\right)^{p/2}\\
\aad \ \le KE\left( \int^{n}_{n-1}
|X(s^-)|^2ds\right)^{p/2}\\
\aad \ \le K\int^n_{n-1}E|X(s^-)|^pds.\eea Given that $X(t)$ is
exponential $p$-stable, we have
\[E|X(t)|^p\le K|x|^pe^{-\kappa t}\text{ for all }  t\ge 0.\]
Substituting the above bounds to \eqref{4.3}, careful calculations
lead to
\[E\left[\sup_{n-1\le t\le n}
|X(t)|^p\right]\le K e^{-\kappa(n-1)}.\]

For any $1\le p<2$, we have
\[E\left[\sup_{n-1\le t\le n}|X(t)|^p \right]\le
\left(E (\sup_{n-1\le t\le n}|X(t)|^{2p}) \right)^{1/2} \le
Ke^{-\kappa (n-1)}.\] Finally, for any $0<p<1$, we have
\[|X(t)|^p=|X(t)|^pI_{\{X(t)\ge 1\}}
+|X(t)|^pI_{\{X(t)<1\}}\le 1+|X(t)|^{1+p}.\] Therefore,
\[E\left(\sup_{n-1\le t\le n}|X(t)|^p\right)
\le 1+E(\sup_{n-1\le t\le n} |X(t)|^{1+p})\le Ke^{-\kappa(n-1)}.\]
Note that in the above $K$ and $\kappa$ have different values in
different appearances. Now let $\varepsilon\in(0,\kappa)$ be
arbitrary, then it follows that \bea \ad P\left\{\sup_{n-1\le t\le
n}|X(t)|^p>e^{-(\kappa-\varepsilon)(n-1)}\right\}\\
\aad \ \le e^{(\kappa-\varepsilon)(n-1)}E\left(\sup_{n-1\le
t\le n}|X(t)|^p\right)\\
\aad \ \le Ke^{-\varepsilon(n-1)} .\eea Since
$\sum^\infty\limits_{n=1} K\exp(-\e (n-1)) < \infty,$
 by Borel-Cantelli lemma, we have
\begin{equation}\label{4.5}
  \sup_{n-1\le t\le n}
  |X(t)|^p\le e^{-(\kappa-\varepsilon)(n-1)} \ \hbox{ a.s.}
\end{equation}
 for all but finitely many $n$.
 Hence, there exists such an
$n_0$ that whenever $n\ge n_0$, \eqref{4.5} holds a.s. So,
\beq{4.6}\frac{1}{t}\ln|X(t)|=\frac{1}{pt} \ln(|X(t)|^p)\le-\frac
{(\kappa-\varepsilon)(n-1)}{p(n-1)} <0\ \hbox{ a.s.}\eeq
Taking $\limsup$ in \eqref{4.6} leads to
almost surely exponential
stability.
Thus, the proof is completed. \qed

\section{Smooth-Dependence  on Initial Data}\label{sec:Lya}
One of the important properties of a diffusion processes is the
continuous and smooth dependence on the initial data. This property
is preserved for the switching diffusion processes with
state-dependent switching; however much work is needed. In what
follows, we show that this property is also preserved for the
switching jump diffusion processes. The results are stated for
multi-dimensional cases, whereas the proofs are carried out for a
one-dimensional process for notational simplicity. Let
$(X(t),\alpha(t))$ denote the switching jump process with initial
condition $(x,\alpha)$ and $(\wdt{X}(t),\wdt{\alpha}(t))$ be the
process starting from $(\wdt x, \alpha)$, let $\Delta\neq 0$ be
small and denote $\tilde x =x+\Delta$ in the sequel.

\begin{lem}\label{lema5.2}
Under conditions \rm{(A1)}-\rm{(A3)}, we have for $0\le t \le
T$ and any positive constant $\iota$, $E|X(t)|^\iota\le |x|^\iota
e^{\kappa t}\le C$,  for $x\neq 0$, $\alpha \in \mathcal{M},$ where
$\kappa=\kappa(\iota,K_0,m,g^*(i))$ and $C=C(\kappa,T)$.
\end{lem}

\para{Proof.} For each $i\in \mathcal{M}$ and $x\neq 0$,
define $V(x,i)=|x|^\iota$ for any $\iota \in \mathbb{R}_+-\{0\}$.
Then for any $\Delta >0$ and $|x|>\Delta$,
 \bea \mathcal{G}|x|^\iota=
\ad \iota|x|^{\iota-2}x'b(x,i)+ \lambda\int_\Gamma(|x+g(x,i,\gamma)|^\iota-|x|^\iota)\pi(d\gamma)\\
\ad+\frac{1}{2}\tr[\sigma(x,i)\sigma'(x,i)\iota|x|^{\iota-4}(|x|^2I+(\iota-2)xx')].\\
 \eea Since $0$ is an equilibrium point,
Cauchy-Schwartz inequality implies $|x'b(x,i)|\le |x||b(x,i)|\le
K_0|x|^2$, \bea \ad \tr(\sigma\sigma')=|\sigma|^2\le K_0|x|^2,
\\ \ad \tr(\sigma
\sigma'xx')=x'\sigma \sigma'x \le |x|^2|\sigma|^2\le K_0|x|^4. \eea
 Therefore, we have \bea  |\mathcal{G}|x|^\iota|\le
\ad K_0\iota|x|^\iota+\frac{1}{2}K_0\iota|x|^{\iota-2}(|x|^2+(\iota-2)|x|^2)\\
\ad+ \lambda|x|^\iota(|1+g^*(i)|^\iota-1)\le \kappa|x|^\iota. \eea
Define the stopping time $ \tau_{\Delta}:=\inf\{t\ge 0, |X(t)|\le
\Delta\}.$ Then by the generalized It\^{o} lemma, we obtain \bea
E|X(\tau_{\Delta}\wedge
t)|^\iota\ad=|x|^\iota+E\int^{\tau_{\Delta}\wedge t}_0
\mathcal{G}|X(u)|^\iota du\\
\ad \le|x|^\iota+\kappa E\int^{\tau_{\Delta}\wedge t}_0 |X(u)|^\iota
du\\
\ad \le|x|^\iota+\kappa E\int^t_0 |X(u\wedge\tau_\Delta)|^\iota du.
\eea

By Gronwall's inequality, it follows that \bea
E|X(\tau_{\Delta}\wedge t)|^\iota\le |x|^\iota e^{\kappa t}.\eea
Letting $\Delta\to 0$, by virtue of non-zero property of $X(t)$
shown in \lemref{0-inacc}, we have \bea E|X(t)|^\iota \le|x|^\iota
e^{\kappa t}.\eea For $0\le t\le T$, we further have \bea
E|X(t)|^\iota \le|x|^\iota e^{\kappa t}\le |x|^\iota e^{\kappa
T}=C.\eea Thus, the proof is completed. \qed

\begin{thm}\label{eqq}
Under  the conditions of \lemref{lema5.2}, define \beq{5.1} \barray
\disp \phi^{\Delta }(t)\ad=\frac{1}{\Delta }\int^t_0[b(\wdt
{X}(s),\wdt
{\alpha}(s))-b(\wdt  {X}(s),\alpha(s))]ds\\
\aad \ +\frac{1}{\Delta }\int^t_0[\sigma(\wdt {X}(s),\wdt
{\alpha}(s))-\sigma(\wdt  {X}(s),\alpha(s))]dw(s) \\
\aad \ +\frac{1}{\Delta }\int^t_0\int_\Gamma[g(\wdt {X}(s^-),\wdt
{\alpha}(s^-),\gamma)-g(\wdt
{X}(s^-),\alpha(s^-),\gamma)]N(ds,d\gamma). \earray \eeq
 Then we have
$ \lim\limits_{\Delta  \to 0} E\sup\limits_{0\le t \le
T}|\phi^\Delta (t)|^2=0.$
\end{thm}

\para{Proof.} It can be verified that
\beq{5.2} \barray
  E\sup\limits_{0\le t\le T}|\phi^\Delta (t)|^2\ad=
\disp  \frac{K}{\Delta ^2}E\int^T_0|b(\wdt {X}(s),\wdt {\alpha}(s))
  -b(\wdt {X}(s),\alpha(s))|^2ds\\
  \aad \ +\frac{K}{\Delta ^2}E\sup\limits_{0\le t\le T}|\int^t_0[\sigma(\wdt {X}(s),\wdt {\alpha}(s))
  -\sigma(\wdt {X}(s),\alpha(s))]dw(s)|^2\\
  \aad \ +\frac{K}{\Delta ^2}E\int^T_0\int_\Gamma|g(\wdt {X}(s^-
),\wdt {\alpha}(s^-),\gamma)
  -g(\wdt {X}(s^-),\alpha(s^-),\gamma)|^2ds \pi(d\gamma)\\
  \aad \ +\frac{K}{\Delta ^2}E\sup\limits_{0\le t\le T}|\int^t_0\int_\Gamma[g(\wdt {X}(s^-
),\wdt {\alpha}(s^-),\gamma)
  -g(\wdt {X}(s^-),\alpha(s^-),\gamma)]\wdt {N}(ds,d\gamma)|^2.
\earray \eeq
 Let us first consider the next to the last line of \eqref{5.2}. By
choosing $\eta=\Delta ^{\gamma_0}$ with $\gamma_0>2$ and partition
the interval $[0,T]$ by $\eta$ we obtain
\begin{equation}\label{eq5.3}
\begin{array}{ll}
\ad E\int^T_0\int_\Gamma |g(\wdt {X}(s^-),\wdt
{\alpha}(s^-),\gamma)-g(\wdt
{X}(s^-),\alpha(s^-),\gamma)|^2ds\pi(d\gamma)\\
\aad \ =E \sum^{\lfloor
{\frac{T}{\eta}\rfloor}-1}_{k=0}\int^{k\eta+\eta}_{k\eta}\int_\Gamma
|g(\wdt {X}(s^-),\wdt  {\alpha}(s^-),\gamma)-g(\wdt
{X}(s^-),\alpha(s^-),\gamma)|^2ds\pi(d\gamma)\\
\aad \ =KE \sum^{\lfloor
{\frac{T}{\eta}\rfloor}-1}_{k=0}[\int^{k\eta+\eta}_{k\eta}\int_\Gamma
|g(\wdt {X}(s^-),\wdt  {\alpha}(s^-),\gamma)-g(\wdt
{X}(k\eta),\wdt \alpha(s^-),\gamma)|^2 ds\pi(d\gamma)\\
\aad \hspace*{0.8in} +\int^{k\eta+\eta}_{k\eta}\int_\Gamma
|g(\wdt {X}(k\eta),\wdt  {\alpha}(s^-),\gamma)-g(\wdt
{X}(k\eta),\alpha(s^-),\gamma)|^2ds\pi(d\gamma)\\
\aad \hspace*{0.8in} +\int^{k\eta+\eta}_{k\eta}\int_\Gamma |g(\wdt
{X}(k\eta),{\alpha}(s^-),\gamma)-g(\wdt
{X}(s^-),\alpha(s^-),\gamma)|^2ds\pi(d\gamma)].
\end{array}
\end{equation}
For the third line of \eqref{eq5.3}, we have the following bound by
virtue of \rm{(A2)} and  \cite[Theorem 3.7.1]{Kushner},
\begin{equation}
\begin{array}{ll}
\ad E\int^{k\eta+\eta}_{k\eta}\int_\Gamma
|g(\wdt {X}(s^-),\wdt {\alpha}(s^-),\gamma)-g(\wdt
{X}(k\eta),\wdt  \alpha(s^-),\gamma)|^2ds\pi(d\gamma)\\
\aad \ \le K\int^{k\eta+\eta}_{k\eta}E\left|\wdt {X}(s^-)-\wdt
{X}(k\eta) \right|^2ds \\
\aad \ \le K\int^{k\eta+\eta}_{k\eta}(s-k\eta)ds\le K\eta^2.
\end{array}
\end{equation}
Recall that $K$ is a generic positive constant, whose values may be
different for different appearances. We can derive the upper bound
for the last line of \eqref{eq5.3} similarly,  \bea \ad
E\int^{k\eta+\eta}_{k\eta}\int_\Gamma |g(\wdt
{X}(k\eta),{\alpha}(s^-),\gamma)-g(\wdt
{X}(s^-),\alpha(s^-),\gamma)|^2ds\pi(d\gamma)\le O(\eta^2). \eea

 To treat the term on the next to the last line of \eqref{eq5.3},
 note  that
\beq {eq5.4} \barray \ad
E\int^{k\eta+\eta}_{k\eta}\int_\Gamma |g(\wdt {X}(k\eta),\wdt
{\alpha}(s^-),\gamma)-g(\wdt
{X}(k\eta),\alpha(s^-),\gamma)|^2ds\pi(d\gamma)\\
\aad \ \le KE\int^{k\eta+\eta}_{k\eta}\int_\Gamma |g(\wdt
{X}(k\eta),\wdt  {\alpha}(s^-),\gamma)-g(\wdt
{X}(k\eta),\wdt \alpha(k\eta),\gamma)|^2ds\pi(d\gamma)\\
\aad \quad +KE\int^{k\eta+\eta}_{k\eta}\int_\Gamma |g(\wdt
{X}(k\eta),\wdt  {\alpha}(k\eta),\gamma)-g(\wdt
{X}(k\eta),\alpha(s^-),\gamma)|^2ds\pi(d\gamma). \earray \eeq
 For the term on the second line of \eqref{eq5.4} and
 $k=0,1,\cdots,\lfloor \frac{T}{\eta}\rfloor -1$,
\bea\disp \ad E\int^{k\eta+\eta}_{k\eta}\int_\Gamma |g(\wdt
{X}(k\eta),\wdt {\alpha}(s^-),\gamma)-g(\wdt {X}(k\eta),\wdt
\alpha(k\eta),\gamma)|^2 ds\pi(d\gamma)\\ \aad \
=E\int^{k\eta+\eta}_{k\eta}\int_\Gamma |g(\wdt {X}(k\eta),\wdt
{\alpha}(s^-),\gamma)-g(\wdt {X}(k\eta),\wdt
\alpha(k\eta),\gamma)|^2I_{\{\wdt  \alpha(s^-)\neq \wdt
\alpha(k\eta)\}}ds \pi (d\gamma) \eea
\beq{eq5.6}\barray
 \ad\ \ =E\sum_{i\in
\mathcal{M}}\sum_{j\neq i}\int^{k\eta+\eta}_{k\eta}\int_\Gamma
|g(\wdt {X}(k\eta),j,\gamma)-g(\wdt {X}(k\eta),i,\gamma)|^2I_{\{\wdt
\alpha(s^-)=\wdt \alpha(s)=j\}}I_{\{\wdt \alpha(k\eta)=i\}}ds \pi
(d\gamma)
\\
\aad \ \le KE\sum_{i\in \mathcal{M}}\sum_{j\neq
i}\int^{k\eta+\eta}_{k\eta} [1+|\wdt  {X}(k\eta)|^2]I_{\{\wdt
{\alpha}(k\eta)=i\}}\times E[I_{\{ \wdt  {\alpha}(s)=j\}}|\wdt
{X}(k\eta),\wdt  {\alpha}(k\eta)=i]ds
\\
\aad \ \le KE\sum_{i\in \mathcal{M}}\int^{k\eta+\eta}_{k\eta} [1
+|\wdt {X}(k\eta)|^2]I_{\{\wdt  {\alpha}(k\eta)=i\}}\times
[\sum_{j\neq i} q_{ij}(\wdt {X}(k\eta))(s-k\eta)+o(s-k\eta)]ds
\\
\aad \ \le K\int^{k\eta+\eta}_{k\eta}O(\eta)ds\le K\eta^2.
\earray
\eeq

 In the above, we
employed the fact that the time of jump of $X(t)$  does not coincide
with that of switching part $\alpha(t)$ in \cite[Proposition
2.2,]{Xi}. Also, \lemref{lema5.2} and boundedness
 of $Q(x)$ are involved.
Now let us deal with the last line of \eqref{eq5.4} by using the
basic coupling techniques \cite[p. 11]{Chen04}. Consider the measure
\bea \Lambda((x,j),(\wdt {x},i))=|x-\wdt {x}|+d(j,i),\ \hbox{ where
}\
 d(j,i) = \left\{ \begin{array}{l}
 0\ \hbox{ if } \ j = i, \\
 1\ \hbox{ if } \ j \ne i.
 \end{array} \right.\eea
Let $(\alpha(t),\wdt  \alpha(t))$ be a random process with
a finite state space $\M\times\M$ such that \bea
 \ad P[(\alpha (t + h),\wdt  \alpha (t + h)) = (j,i)|(\alpha (t),\wdt
 \alpha (t)) = (k,l),(X(t),\wdt  X(t)) = (x,\wdt  x)] \\
 \ad = \left\{ \barray
 \wdt  q_{(k,l)(j,i)} (x,\wdt  x)h + o(h),\     \ad\text{ if } \ (k,l) \ne (j,i), \\
 1 + \wdt  q_{(k,l)(k,l)} (x,\wdt  x)h + o(h),\ \ad\text{ if } \
 (k,l) = (j,i), \\
 \earray \right.
 \eea
where $h \to 0$, and the matrix $(\wdt {q}_{(k,l)(j,i)}(x,\wdt
{x}))$ is the basic coupling of matrices $Q(x)=(q_{kl}(x))$ and
$Q(\wdt {x})=(q_{kl}(\wdt {x}))$ satisfying \beq{eq5.7} \barray \wdt
{Q}(x,\wdt {x})\wdt {f}(k,l)\ad=\sum_{(j,i)\in \mathcal{M}\times
  \mathcal{M}}\tilde q_{(k,l)(j,i)}(x,\wdt {x})(\wdt {f}(j,i)-\wdt {f}(k,l))\\
  \ad=\sum
  _j(q_{kj}(x)-q_{lj}(\wdt {x}))^+(\wdt {f}(j,l)-\wdt {f}(k,l))
  \\\aad\ +\sum
  _j(q_{lj}(\wdt {x})-q_{kj}(x))^+(\wdt {f}(k,j)-\wdt {f}(k,l))
  \\\aad\ +\sum
  _j(q_{kj}(x)\wedge
  q_{lj}(\wdt {x}))(\wdt {f}(j,j)-\wdt {f}(k,l))
\earray \eeq for any function $\wdt {f}(\cdot,\cdot)$ defined on
$\mathcal{M}\times \mathcal{M}$. Then we have
\beq{eq5.8}
  \barray
 \ad E[I_{\{\alpha(s)=j\}}|\alpha(k\eta)=i_1,\wdt {\alpha}(k\eta)=i,X(k\eta)=x,\wdt
  X(k\eta)=\wdt {x}]
  \\\ad=\sum_{l\in
\mathcal{M}}E[I_{\{\alpha(s)=j\}}I_{\{\wdt
{\alpha}(s)=l\}}|\alpha(k\eta)=
  i_1,\wdt {\alpha}(k\eta)=i,X(k\eta)=x,\wdt {X}(k\eta)=\wdt {x}]
  \\\ad=\sum_{l\in
  \mathcal{M}}\wdt {q}_{(i_1,i)(j,l)}(x,\wdt {x})(s-k\eta)+o(s-k\eta)=O(\eta).
  \earray
\eeq
 Therefore, for $ k=1,\cdots,\lfloor
\frac{T}{\eta}\rfloor -1$, we have \beq{eq5.9}
  \barray
 \ad E\int^{k\eta+\eta}_{k\eta}\int_\Gamma
|g(\wdt {X}(k\eta),\wdt  {\alpha}(k\eta),\gamma)-g(\wdt
{X}(k\eta),\alpha(s^-),\gamma)|^2ds\pi(d\gamma)
\\
\aad \ =E\sum_{i\in \mathcal{M}}\sum_{j\neq
i}\int^{k\eta+\eta}_{k\eta}\int_\Gamma |g(\wdt {X}(k\eta),i,\gamma)
-g(\wdt  {X}(k\eta),j,\gamma)|^2
I_{\{\alpha(s)=\alpha(s^-)=j\}}I_{\{\wdt
{\alpha}(k\eta)=i\}}ds\pi(d\gamma)
\\
\aad \ \le KE\sum_{i,i_1\in \mathcal{M}}\sum_{j\neq
i}\int^{k\eta+\eta}_{k\eta}[1+|\wdt {X}(k\eta)|^2]I_{\{\wdt {\alpha}(k\eta)=i,
\alpha(k\eta)=i_1\}}
\\
\aad \quad \times E[I_{\{\alpha(s)=j\}}|\alpha(k\eta)=i_1,\wdt
{\alpha}(k\eta)=i,X(k\eta)=x,\wdt {X}(k\eta)=\wdt {x}]ds=O(\eta^2).
  \earray
\eeq For $k=0$, recall that $\alpha(0)=\wdt{\alpha}(0)=\alpha$,
$X(0)=x$ and $\wdt {X}(0)=\wdt {x}$, we have \beq{eq510}
  \barray
   \ad E\int^\eta_0\int_\Gamma|g(\wdt {X}(0),\wdt {\alpha}(0),
    \gamma)-g(\wdt {X}(0),\alpha(s),\gamma)|^2ds\pi(d\gamma)
     \\
     \aad \ =E\int^\eta_0\int_\Gamma\sum_{j\neq
    \alpha}|g(\wdt {x},\alpha,\gamma)-
g(\wdt {x},j,\gamma)|^2I_{\{\alpha(s)=j\}}ds\pi(d\gamma)
    \\
    \aad \ \le K\sum_{j\neq \alpha}\int^\eta_0 [1+\tilde{x}^2]E[I_{\{\alpha(s)=j\}}|\alpha(0)=\alpha,\wdt{X}(0)=\wdt x]ds\\
    \aad \ \le K\int^\eta_0\sum_{j\neq
    \alpha}[q_{\alpha
    j}(\wdt {x})s+o(s)]ds\le K\eta^2.
  \earray
\eeq Thus, for $ k=0,1,\cdots,\lfloor \frac{T}{\eta}\rfloor -1$,
\beq{eq5.11} E\int^{k\eta+\eta}_{k\eta}\int_\Gamma |g(\wdt
{X}(k\eta),\wdt  {\alpha}(k\eta),\gamma)-g(\wdt
{X}(k\eta),\alpha(s^-),\gamma)|^2ds\pi(d\gamma)\le K\eta^2. \eeq

Now we can obtain
 \bea \ad E\int^T_0\int_\Gamma |g(\wdt {X}(s^-),\wdt
{\alpha}(s^-),\gamma)-g(\wdt {X}(s^-),\alpha(s^- ),\gamma)|^2ds
\pi(d\gamma) \le\sum^{\lfloor{\frac{T}{\eta}}\rfloor-1}_{k=0}
K\eta^2\le K\eta. \eea Likewise, we also obtain the bound for the
martingale part
 \bea \ad
 E\sup\limits_{0\le t\le T}|\int^t_0\int_\Gamma[g(\wdt {X}(s^-
),\wdt {\alpha}(s^-),\gamma)-g(\wdt{X}(s^-),\alpha(s^-),\gamma)]\wdt
{N}(ds,d\gamma)|^2 \le K\eta.\eea
For the drift and diffusion parts
involved, the argument in \cite[Lemma 4.3]{YinZhu} leads to \bea \ad
E\int^T_0|(b(\wdt {X}(s),\wdt {\alpha}(s))-b(\wdt
{X}(s),\alpha(s))|^2ds
\le K\eta, \\
\ad E\sup\limits_{0\le t\le T}|\int^t_0[\sigma(\wdt {X}(s),\wdt
{\alpha}(s)) -\sigma(\wdt {X}(s),\alpha(s))]dw(s)|^2\le K\eta.\eea
Therefore, we obtain \beq{5.121}
 E\sup_{0\le t\le
T}|\phi^\Delta (t)|^2\le K\frac{\eta}{\Delta
^2}=K\Delta^{\gamma_0-2}\to 0 \text{ as }\Delta \to 0.\eeq
 This concludes the proof. \qed

\begin{lem}\label{lem5.22}
  Under the conditions of \thmref{eqq},
  $ E[\sup\limits_{0\le t \le
  T}|\wdt {X}^{\wdt {x},\alpha}(t)-X^{x,\alpha}(t)|^2]\le
  C|\wdt {x}-x|^2,$
   where the constant $C$ satisfies $C=C(K_0,T)$.
\end{lem}

\para{Proof.} Let $T>0$ be fixed and recall that $\Delta =\wdt {x}-x$, then we have
$\wdt {X}^{\wdt x, \alpha}(t)-X^{x,\alpha}(t)=\Delta +A(t)+B(t),$
where \beq{eq5.12}
  \barray
  A(t)\ad =\int^t_0[b(\wdt {X}(s),\wdt {\alpha}(s))-
b(\wdt {X}(s),\alpha(s))]ds\\
  \aad \ +\int^t_0[\sigma(\wdt {X}(s),\wdt {\alpha}(s))-
\sigma(\wdt {X}(s),\alpha(s))]dw(s)\\
  \aad \ +\int^t_0\int_\Gamma[g(\wdt {X}(s^-),\wdt {\alpha}(s^-),
  \gamma)-g(\wdt {X}(s^-),\alpha(s^-),\gamma)]N(ds,d\gamma)\\
  \ad =\Delta \phi^\Delta (t),
  \earray
\eeq
\beq{eq5.13}
  \barray
  B(t)\ad=\int^t_0[b(\wdt {X}(s),\alpha(s))-b(X(s),\alpha(s))]ds\\
  \aad \ +\int^t_0[\sigma(\wdt {X}(s),\alpha(s))-\sigma(X(s),\alpha(s))]dw(s)\\
  \aad \ +\int^t_0\int_\Gamma[g(\wdt {X}(s^-),\alpha(s^-),\gamma)-g(X(s^-),\alpha(s^-
),\gamma)]N(ds,d\gamma). \earray \eeq
 Hence
 \bea
\ad \sup_{t\in[0,T]}|\wdt {X}^{\wdt x, \alpha}(t)-X^{x,\alpha}(t)|^2
\le 3\Delta ^2+3\sup_{t\in[0,T]}|A(t)|^2+3\sup_{t\in[0,T]}
|B(t)|^2.\eea

 It follows from \eqref{5.121}
that \bea E[\sup\limits_{t\in[0,T]}|A(t)|^2]\le \Delta^2
E[\sup\limits_{t\in[0,T]}|\phi^\Delta(t)|^2]\le K
\Delta^{r_0}=o(\Delta^2).\eea By the H\"{o}lder inequality and the
Lipschitz continuity, we have \bea \ad
E[\sup_{t\in[0,T]}|\int^t_0[b(\wdt
{X}(s),\alpha(s))-b(X(s),\alpha(s))]ds|^2] \le K\int^T_0E|\wdt
{X}(s)-X(s)|^2ds\eea and
 \bea
\ad E[\sup_{t\in[0,T]}|\int^t_0\int_\Gamma[
 g(\wdt {X}(s^-),\alpha(s^-),\gamma)-g(X(s^-),\alpha(s^-
),\gamma)]ds\pi(d\gamma)|^2]\\
 \aad \
\le K\int^T_0E|\wdt {X}(s^-)-X(s^-)|^2ds.\eea
 Then the basic
properties of stochastic integrals (w.r.t. $w\cd$ and $\wdt{N}\cd$)
together with the Lipschitz continuity lead to \bea\ad
E[\sup_{t\in[0,T]}|\int^t_0[\sigma(\wdt
{X}(s),\alpha(s))-\sigma(X(s),\alpha(s))]dw(s)|^2] \le
K\int^T_0E|\wdt {X}(s)-X(s)|^2ds\eea
 and
 \bea\ad
 E[\sup_{t\in[0,T]}|\int^t_0\int_\Gamma [g(\wdt {X}(s^-),\alpha(s^-),\gamma)-g(X(s^-),
 \alpha(s^-),\gamma)]\wdt {N}(ds,d\gamma)|^2]\\
 \aad \
\le K\int^T_0E|\wdt {X}(s^-)-X(s^-)|^2ds.\eea So, \beq{eq5.14}
\barray \ad E[\sup_{t\in[0,T]
  }|\wdt {X}^{\wdt {x},\alpha}(t)-X^{x,\alpha}(t)|^2]
   \le 3\Delta ^2+K\int^T_0E[\sup_{u\in[0,T]}|\wdt {X}(u)-X(u)|^2]du+o(\Delta^2).
\earray \eeq

Now, by Gronwall's inequality \bea\ad
E[\sup_{t\in[0,T]}|\wdt{X}^{\wdt
{x},\alpha}(t)-X^{x,\alpha}(t)|^2]\le
  3\Delta ^2\exp(KT)+o(\Delta^2)\le K|\wdt {x}-x|^2.\eea
  Thus, we have completed the proof. \qed

Let us introduce some notations to proceed. Recall that a vector
 $\beta=(\beta_1,\beta_2,\cdots, \beta_r)$ with nonnegative integer component is referred to as a multi-index.
 Put
$ |\beta|=\beta_1+\beta_2+\cdots+\beta_r,$
  we define $D^\beta_x$ as $$
 D^\beta_x=\frac{\partial^\beta}{\partial
x^\beta}=\frac{\partial^{|\beta|}}{\partial x^{\beta_1}_1 \cdots
\partial^{\beta_r}_{x_r}}. $$
Recall that $\Delta=\wdt
x-x$ and define \beq{eq5.161}
 Z^\Delta (t)=
\frac{\wdt {X}^{\wdt x,\alpha}(t)-X^{x,\alpha}(t)}{\Delta }. \eeq
Then we have the following expression: \beq{eq5.16} \barray
 \disp Z^\Delta(t) \ad  =1+\phi^{\Delta }(t)+
\frac{1}{\Delta }\int^t_0[b(\wdt {X}(s),\alpha(s))-b(X(s),\alpha(s))]ds\\
\aad\quad + \frac{1}{\Delta }\int^t_0[\sigma(\wdt {X}(s),
\alpha(s))-\sigma(X(s),\alpha(s))]dw(s)\\
\aad\quad +\frac{1}{\Delta }\int^t_0\int_\Gamma[g(\wdt
{X}(s^-),\alpha(s^-),\gamma)-
g(X(s^-),\alpha(s^-),\gamma)]N(ds,d\gamma) ,\earray \eeq
 where $\phi^\Delta(t)$ is defined in \eqref{5.1}.

\begin{lem}\label{eq}
Under the conditions of \thmref{lem5.22}, assume that for each $i\in
\M$, $b(\cdot, i)$, $\sigma(\cdot,i)$ and $g(,i,\gamma)$ have
continuously partial derivatives with respect to the variable $x$ up
to the second order and that \bea |D^\beta_x b(x,i)|+|D^\beta_x
\sigma(x,i)|+|D^\beta_x g(x,i,\gamma)|\le K(1+|x|^\rho),\eea
 where $K$ and $\rho$ are positive constants and $\beta$ is a multi-index
with $|\beta|\le 2$. Then $X^{x,\alpha}(t)$ is twice continuously
differentiable in mean square with respect to $x$.
\end{lem}

\para{Proof.} Given the
definition of $Z^{\Delta}(t)$ above and \thmref{eqq}, we just need
to consider the last three terms of \eqref{eq5.16}. First, note that
 \beq{eq5.17} \barray
 \ad \frac{1}{\Delta }\int^t_0\int_\Gamma[g(\wdt {X}(s^-),
\alpha(s^-),\gamma)-g(X(s^-),\alpha(s^-),\gamma)]ds\pi(d\gamma)
\\
\aad \ =\frac{1}{\Delta }\int_\Gamma\int^t_0\int^1_0
\frac{d}{{d\nu}}g(X(s^-)+\nu(\wdt {X}(s^-)-X(s^-)),\alpha(s^-),
\gamma)d\nu ds\pi(d\gamma)
\\
\aad \ =\int_\Gamma\int^t_0[\int^1_0 g_x(X(s^-)+\nu(\wdt
{X}(s^-)-X(s^-)),\alpha(s^-),\gamma)d\nu]Z^\Delta
(s^-)ds\pi(d\gamma), \earray \eeq
 where $g_x(\cdot)$ denotes the partial derivative of $g(\cdot,i,\gamma)$ with respect to $x$.  It follows from
\lemref{lem5.22}
 that for any $s \in [0, T]$,
$\wdt {X}(s^-)-X(s^-)\to 0$  in probability as $\Delta
\to 0.$ This implies that
\begin{equation}\label{5.19}
\int^1_0 g_x(X(s^-)+\nu(\wdt
{X}(s^-)-X(s^-)),\alpha(s^-),\gamma)d\nu \to
g_x(X(s^-),\alpha(s^-),\gamma) \end{equation} in probability as
$\Delta  \to 0$. Therefore, we have
\beq{eq5.20}\barray\disp\frac{1}{\Delta}\int^t_0\int_\Gamma[g(\wdt{X}(s^-),\alpha(s^-),\gamma)-g(X(s^-),\alpha(s^-),\gamma)]ds\pi
(d\gamma)\\\disp\rightarrow\int_0^t\int_\Gamma
g_x(X(s^-),\alpha(s^-),\gamma)Z^\Delta(s^-)ds\pi
(d\gamma).\earray\eeq Similarly, we have
 \beq {b}
\frac{1}{\Delta}\int^t_0
[b(\wdt{X}(s),\alpha(s))-b(X(s),\alpha(s))]ds\to \int^t_0
b_x(X(s),\alpha(s))Z^\Delta(s)ds \eeq in probability as $\Delta\to
0$  and \beq{sigma}
 \frac{1}{\Delta}\int^t_0
[\sigma(\wdt{X}(s),\alpha(s))-\sigma(X(s),\alpha(s))]dw(s)\to
\int^t_0\sigma_x(X(s),\alpha(s))Z^\Delta(s)dw(s) \eeq in probability
as $\Delta \to 0$. $b_x(\cdot)$ and $\sigma_x(\cdot)$ denote the
partial derivative of $b(\cdot, i)$ and $\sigma(\cdot, i)$ with
respect to $x$, respectively.
Recall the definition of $Z^\Delta(t)$
in equation \eqref{eq5.161}, \thmref{eqq},
\eqref{5.19}-\eqref{sigma}, and
 \cite[Theorem 5.5.2]{F} yield
 \beq{zt} E|Z^\Delta (t)-\varsigma(t)|^2 \to 0
 \text{ as }
 \Delta  \to 0.\eeq
 where
 \beq{eq5.21}
 \barray
 \varsigma(t)\ad=1+\int^t_0 b_x(X(s),\alpha(s))\varsigma(s)ds+\int^t_0
 \sigma_x(X(s),\alpha(s))\varsigma(s)dw(s)
 \\
 \aad \ +\int^t_0\int_\Gamma g_x(X(s^-),\alpha(s^-),\gamma)\varsigma(s^-)N(ds,d\gamma)
 \earray
\eeq
 and $\varsigma(t)=\varsigma^{x,\alpha}(t)$ is mean square continuous with respect to
$x$. Therefore, $\frac{\partial}{\partial x}X^{x,\alpha}(t)$ exists
in the mean square sense and $\varsigma(t)=\frac{\partial}{\partial
x}X^{x,\alpha}(t)$. Likewise, we can show
$\frac{\partial^2}{\partial x^2}X^{x,\alpha}(t)$ exists in the mean
square sense and is mean square continuous with respect to $x$. \qed

\begin{lem}\label{z111}
Under the assumptions of \lemref{eq}, we have $\sup\limits_{t\in
[0,T]}E|\varsigma(t)|^2\le K=K(x, \wdt x, T, K_0)<\infty.$
\end{lem}

\para{Proof.}
For any $t \in [0,T]$, $ E|\varsigma(t)|^2\le
2E|\varsigma(t)-Z^\Delta (t)|^2+2E|Z^\Delta (t)|^2$. By \eqref{zt},
it suffices to consider the last term above. In fact,
\bea E|Z^\Delta (t)|^2 \ad \le K+5E|\phi^\Delta (t)|^2+5E|\frac{1}{\Delta
}\int^t_0[b(\wdt
{X}(u),\alpha(u))-b(X(u),\alpha(u))]du|^2\\
\aad \ +5E|\frac{1}{\Delta }\int^t_0[\sigma(\wdt
{X}(u),\alpha(u))-\sigma(X(u),\alpha(u))]dw(u)|^2\\
\aad \ +5E|\frac{1}{\Delta }\int^t_0\int_\Gamma[g(\wdt
{X}(u^-),\alpha(u^-),\gamma)-g(X(u^-),
\alpha(u^-),\gamma)]N(du,d\gamma)|^2,\eea
so
 \beq{eq5.22}
  \barray
E|Z^\Delta (t)|^2
\ad  \le K+5t\frac{1}{|\Delta |^2}E\int^t_0|b(\wdt
{X}(u),\alpha(u))-b(X(u),\alpha(u))|^2du\\
 \aad \ +5\frac{1}{|\Delta |^2}E\int^t_0|\sigma(\wdt
{X}(u),\alpha(u))-\sigma(X(u),\alpha(u))|^2du\\
\aad \ +5t\frac{1}{|\Delta |^2}E\int^t_0\int_\Gamma|g(\wdt
{X}(u^-),\alpha(u^-),\gamma)-g(X(u^-),\alpha(u^-),\gamma)|^2du\pi(d\gamma)+\\
\aad \ +5\frac{1}{|\Delta |^2}E\int^t_0\int_\Gamma|g(\wdt
{X}(u^-),\alpha(u^-),\gamma)-g(X(u^-),\alpha(u^-),\gamma)|^2du\pi(d\gamma)\\
\ad \le
K+5K_0(T+1)\frac{1}{|\Delta |^2}E\int^t_0|\wdt {X}(u)-X(u)|^2du\\
\aad \ +5K_0(T+1)\frac{1}{|\Delta |^2}E\int^t_0|\wdt {X}(u^-)-X(u^-)|^2
du \le K=K(x,\wdt x,T,K_0).
  \earray
\eeq
 Hence the proof is completed. \qed

\begin{lem}\label{5.7}
Assume the conditions of \lemref{z111} hold. Then the function
$E|X^{x,\alpha}(t)|^p$ is twice continuously differentiable with
respect to the variable $x$, except possibly at $x=0$.
\end{lem}

\para{Proof.} In what follows, let
$u(t,x,\alpha)=E[\phi(X(t),\alpha(t))]=E|X^{x,\alpha}(t)|^p$, then
\bea \disp \frac{u(t,\tilde
x,\alpha)-u(t,x,\alpha)}{\Delta}\ad=\frac{1}{\Delta}E[|\wdt {X}(t)|^p-|X(t)|^p]\\
\ad=\frac{1}{\Delta}E\int^1_0\frac{d}{dv}|X(t)+v(\wdt
{X}(t)-X(t))|^p dv\\
\ad=E[Z^\Delta(t)\int^1_0|X(t)+v(\wdt {X}(t)-X(t))|^p_x dv],\eea
where $|\cdot|^p_x$ denotes the partial derivative of
$\phi(\cdot,i)=|\cdot|^p$ with respect to $x$. Consider \beq{5.27}
 \barray
 \ad |\frac{1}{\Delta}E[|\wdt
{X}(t)|^p-|X(t)|^p]-E[|X(t)|^p_x \varsigma(t)]|\\
\aad \ \le |E\int^1_0[|X(t)+v(\wdt {X}(t)-X(t))|^p_x dv
Z^\Delta(t)]-E|X(t)|^p_x \varsigma(t)|\\
\aad \ \le E\int^1_0|\left[|X(t)+v(\wdt {X}(t)-X(t))|^p_x
dv-|X(t)|^p_x
\right]Z^\Delta(t)|+E||X(t)|^p_x[Z^\Delta(t)-\varsigma(t)]|. \earray
\eeq For the second part of last line of \eqref{5.27}, by
Cauchy-Schwartz inequality, we obtain \bea \disp
E\left||X(t)|^p_x[Z^\Delta(t)-\varsigma(t)]\right| \ad\le
E^{\frac{1}{2}}|X(t)|^{2p}_xE^{\frac{1}{2}}
[Z^\Delta(t)-\varsigma(t)]^2\\
\ad \le KE^{\frac{1}{2}}[Z^\Delta(t)-\varsigma(t)]^2 \to 0 \text{ as
}\Delta\to 0. \eea

 Here we used \lemref{lema5.2} and \eqref{zt}.
Similarly, we can show the first term of last line of
 \eqref{5.27} goes to $0$ as $\Delta \to 0$.
 Thus $E|X^{x,\alpha}(t)|^p$ is differentiable with
 respect to the variable $x$. Likewise,
  we can also see it is twice
  continuously differentiable with
 respect to the variable $x$.
 As a nice application of the smooth
  dependence on the initial data, we
 obtain a Lyapunov converse theorem, namely, necessary conditions
 for exponential $p$ stability.

\begin{thm}\label{th5.1} Assume that the conditions of \lemref{5.7} hold and that the equilibrium point $0$ is
exponentially $p$-stable. Then for each $i \in \mathcal{M}$, there
exists a function $V(\cdot,i)\in C^2(\rr^r: \rr_+)$
such that \bea \ad k_1|x|^p\le V(x,i)\le k_2|x|^p \ x\in D,\\
\ad \mathcal{G}V(x,i)\le -k_3|x|^p \ \text{  for  all  } \ x\in D-\{0\},\\
\ad \l \frac{\partial V}{\partial x_j}(x,i)\r\le k_4|x|^{p-1},\\
\ad \l\frac{\partial^2 V}{\partial x_j
\partial x_l}(x,i)\r \le  k|x|^{p-2}.\eea
for all $1 \le j,l\le r, x\in D-\{0\}$, and for some positive
constants $k, k_1, k_2, k_3$ and $k_4$, where $D$ is a neighborhood
of $0$.
\end{thm}

\para{Proof.} For each $i\in \mathcal{M}$,
consider the function \bea \ad V(x,i)=\int^T_0
E|X^{x,i}(u)|^pdu.\eea It follows from \lemref{5.7}, $V(x,i)$ is
twice continuously differentiable with respect to $x$ except
possibly at $0$. The equilibrium point $0$ is exponential
$p$-stable, therefore there is a $\kappa >0$ such that \bea \ad
V(x,i)=\int^T_0 E|X^{x,i}(u)|^pdu \le \int^T_0 K|x|^p e^{- \kappa
u}du \le k_2|x|^p. \eea  For the function $|x|^p$, we have
$|\mathcal{G}|x|^p|\le K|x|^p$ for some positive real number $K$. An
application of generalized It\^{o}'s formula leads to \bea \disp
E|X^{x,i}(T)|^p-|x|^p=E\int^T_0\mathcal{G}|X^{x,i}(u)|^pdu
\ge-KE\int^T_0|X^{x,i}(u)|^pdu=-K
V(x,i). \eea Again recall that equilibrium point $x=0$ is
exponential $p$-stable, we can choose $T$ such that $
E|X^{x,i}(T)|^p\le \frac{1}{2}|x|^p,$ and therefore, we have
$V(x,i)\ge \frac{|x|^p}{2K}=k_1|x|^p$. Notice that \bea \disp
\mathcal{G}V(x,i)=\int^T_0 \mathcal{G}E|X^{x,i}(u)|^pdu. \eea Let
$u(t,x,i)=E|X^{x,i}(t)|^p$, by the similar argument in step 1 and
step 2 of \cite[Theorem 7.10]{YinZ10}, we obtain \bea
\mathcal{G}V(x,i)
\ad =\int^T_0 \mathcal{G}E|X^{x,i}(u)|^pdu=u(T,x,i)-u(0,x,i)\\
\ad=E|X^{x,i}(T)|^p-E|X^{x,i}(0)|^p=E|X^{x,i}(T)|^p-|x|^p\\
\ad\le-\frac{1}{2}|x|^p=-k_3|x|^p.\eea

Note that \bea \disp \frac{\partial E|X^{x,i}(t)|^p}{\partial
x_j}=pE |X^{x,i}(t)|^{p-1}\hbox{sgn}(X^{x,i}(t))\frac{\partial
X^{x,i}(t)}{\partial x_j}, \eea so \bea\disp \l\frac{\partial
E|X^{x,i}(t)|^p}{\partial
x_j}\r \ad=pE\(|X^{x,i}(t)|^{p-1}\l\frac{\partial
X^{x,i}(t)}{\partial x_j}\r\) \\
\ad \le p
E^{\frac{1}{2}}|X^{x,i}(t)|^{2p-2}E^{\frac{1}{2}}\l\frac{\partial
X^{x,i}(t)}{\partial x_j}\r^2\\
\ad \le K(|x|^{2p-2}e^{-\kappa
t})^{\frac{1}{2}}=K|x|^{p-1}e^{-\kappa t/2}. \eea For the last line
above, we used the \lemref{lema5.2} and \lemref{z111}. Consequently,
we have
 \bea\disp  \left|\frac{\partial
V(x,i)}{\partial x_j}\right|= \l\int^T_0\frac{\partial}{\partial
x_j}E|X^{x,i}(u)|^p du \r\le\int^T_0 K|x|^{p-1}e^{-\kappa u/2}du \le
k_4|x|^{p-1}. \eea We can have estimate of the second derivative of
$V(x,i)$ by similar argument, the theorem is thus proved. \qed

\section{Asymptotic Stability in Distribution}
\label{sec:as-d}

For practical systems, frequently, we do not have information
regarding the
 equilibria of the systems. Nevertheless, the systems still possesses
certain
kind of stability properties. Thus it is necessary to extend our
definition to consider the so-called the asymptotic stability in
distribution. Here, the assumptions  $b(0,i)=0$ and $\sigma(0,i)=0$
for each item $i\in\M$ are not needed. That is,
 the system under consideration may have no
equilibrium point at all.
To proceed, let us first give two definitions.

\begin{defn}\label{def-std}{\rm
The system given by \eqref{2.1} and \eqref{2.11} is asymptotically
stable in distribution if, there exists such a probability measure
$\nu(\cdot \times \cdot)$ on $\rr^r \times \mathcal {M}$ that the
transition probability $p(t,x,\alpha,dy \times\{i\})$ of
$\left(X(t),\alpha(t)\right)$ converges weakly to $\nu(dy \times
\{i\})$ as $t \rightarrow \infty$ for every $(x,\alpha)\in \rr^r
\times \mathcal {M}$.

}\end{defn}

\begin{defn}\label{def-P1}
{\rm  The definitions of (P1) and (P2) are as follows.

\begin{itemize}
\item
The switching jump diffusion process  given by \eqref{2.1} and
\eqref{2.11} is said to have property ({\bf P1}) if, for any
$(x,\alpha)\in \mathbb{R}^r \times \mathcal {M}$ and any
$\varepsilon > 0$, there exists a constant $R>0$ such that
 \bea
P\{|X^{x,\alpha}(t)| \ge R\} < \varepsilon, \text{ for  any } t\ge
0.
 \eea

\item The switching jump diffusion process
given by \eqref{2.1} and \eqref{2.11} is said to have property ({\bf
P2}) if, for any $\varepsilon
>0$ and any compact subset $\K$ of $\rr^r$,
there exists a $T=T(\varepsilon , \K)>0$ such that \bea
P(|X^{x_0,i_0}(t)-X^{y_0,i_0}(t)|\le \varepsilon)\rightarrow 1
\text{ as } t\rightarrow \infty, \eea whenever $(x_0,y_0,i_0)\in \K
\times \K \times \mathcal {M}$.
\end{itemize}
}\end{defn}

In this section, we first establish asymptotic stability in
distribution of the process $(X(t),\alpha (t))$ in which $\alpha(t)$
is a Markov chain that is independent of the Brownian motion, which
is referred as Markov switching jump diffusions. Then we further
extend the results to state-dependent switching process.

\subsection{Markov Switching Jump Diffusions}
Throughout this section, $\alpha(t)$ is a Markov chain independent
of the Brownian motion.  We first establish a result on stability in
distribution.
\begin{prop}\label{st-in-d}
Suppose that {\rm (A2)} is satisfied,
that $b(\cdot,i)$, $\sg(\cdot,i)$,
and $g(\cdot,i,\ga)$ grow at most linearly for each $i\in \M$ and
$\ga \in \Gamma$, that conditions {\rm (P1)} and {\rm (P2)}
hold, and that the generator of the Markov chain
$Q$ is irreducible. Then the switching jump diffusion process
$(X(t),\al(t))$ is stable in distribution. \end{prop}

\para{Proof.} We note that
 \cite[Theorem 3.1]{Maoyuan}
in fact works not only for Markov switching diffusion processes but
also for more general Markov processes. In our current setup,
$(X(t),\al(t))$ is a Markov process. So we can use essentially the
same steps as in the aforementioned reference to show the process is
stable in distribution. The verbatim argument is omitted. \qed

Our next task is to find sufficient conditions that ensure conditions
(P1) and (P2) are in force. The result is stated in the next theorem.

\begin{thm}
\label{Theorem 6.1} Assume that for each $i\in\M$, there exists
function $V(\cdot,i)
\in C^2(\rr^r: \rr_+)$ satisfying the following two
conditions:
 There exists a positive real
number $\beta$ such that \beq{eq6.1} \mathcal {G} V(x,i) \le -\beta
V(x,i),\eeq
  \beq{cd-2} V_R:=\inf\limits_{|x|\ge R \atop
  i\in \mathcal {M}} V(x,i)\rightarrow \infty
\text{ as } R\rightarrow \infty.\eeq
Then {\rm (P1)} and {\rm (P2)} hold. \end{thm}

\para{Proof.} Let us first verify {\rm (P1)}.
Define the stopping time \bea \tau_R:=\inf\{t\ge 0: |X(t)|\ge
R\}.\eea

Consider $V(x,i)e^{\beta t}$ and let $t_R=\tau_ R\wedge t$. By
virtue of Dynkin's formula, we have \bea
E_{x,\alpha}[V(X(t_R),\alpha(t_R))e^{\beta t_R}]-V(x,\alpha)
  \ad =E_{x,\alpha} \int^{t_R}_0 e^{\beta s}\mathcal
  {G}V(X(s),\alpha(s))ds
  \\ \aad \ +\beta E_{x,\alpha}
  \int ^{t_R}_0 e^{\beta s}V(X(s),\alpha(s))ds,
\eea
where $E_{x,\al}$ denotes the expectation with $X(0)=x$ and
$\al(0)=\al$.

Hence, by virtue of \eqref{eq6.1}, $
 E_{x,\alpha}V(X(t_R),\alpha(t_R))\le
V(x,\alpha)e^{-\beta t_R}.$
We further have
\bea \ad
 V_R P\{\tau_R\le t\}\le
E_{x,\alpha}[V(X(\tau_R), \alpha(\tau_R))I_{\{\tau_R\le t\}}]\le
V(x,\alpha)e^{-\beta \tau_R}. \eea

Note that $\tau_R \le t$ if and only if $\sup\limits_{0 \le u \le
t}|X(u)|\ge R $. Therefore, it follows that
$$P\{\sup_{0 \le u \le t}|X^{x,\alpha}(u)|\ge R\}\le
\frac{V(x,\alpha)e^{-\beta \tau_R}}{V_R} \le \frac
{V(x,\alpha)}{V_R}.$$ Then upon using \eqref{cd-2}, $
 P\{|X^{x,\alpha}(t)|\ge R\}
\to 0$  as  $R \to \infty$, for all
$t \ge 0$. To guarantee {\rm (P2)} hold, similar technique is
involved here.
  But now we need to consider the difference between two solutions of
  equation \eqref{2.1} starting from different initial values in
  compact set $\K$.
  Namely, $(x,\alpha)$ and $(y,\alpha)$.
  \bea\ad X^{x,\alpha}(t)-X^{y,\alpha}(t)\\
  \aad \ =x-y+\int^t_0
  [b(X^{x,\alpha}(s),\alpha(s))-b(X^{y,\alpha}(s),\alpha(s))]ds\\
  \aad\quad +
  \int^t_0[\sigma(X^{x,\alpha}(s),\alpha(s))-\sigma(X^{y,\alpha}(s),\alpha(s))]dw(s)\\
  \aad\quad +
  \int^t_0\int_\Gamma
  [g(X^{x,\alpha}(s^-),\alpha(s^-),\gamma)-g(X^{y,\alpha}(s^-),\alpha(s^-),\gamma)]
  N(ds,d\gamma).
  \eea
Let $Z^{x,y,\alpha}(t)=X^{x,\alpha}(t)-X^{y,\alpha}(t)$, so
$Z(0)=z=x-y$. Then \bea dZ^{x,y,\alpha}(t)\ad =
  [b(X^{x,\alpha}(t),\alpha(t))-b(X^{y,\alpha}(t),\alpha(t))]dt\\
  \aad \ + [\sigma(X^{x,\alpha}(t),
  \alpha(t))-\sigma(X^{y,\alpha}(t),\alpha(t))]dw(t)
  \\
  \aad \ +\int_\Gamma [g(X^{x,\alpha}(t^-),
  \alpha(t^-),\gamma)-g(X^{y,\alpha}(t^-),
  \alpha(t^-),\gamma)] N(dt,d \gamma).
\eea

Define a stopping time $\tau_\varepsilon:=\inf\{t \ge 0, |X^{x,
  \alpha}(t)-X^{y,\alpha}(t)| \ge \varepsilon\}$ and
  let $t_\varepsilon=\tau_\varepsilon \wedge t $.
Then we have \bea \disp E_{z,\alpha}V(Z(t_\varepsilon),
\alpha(t_\varepsilon))-V(z,\alpha)\ad =E_{z,\alpha}
\int^{t_\varepsilon}_0 \mathcal {G}V(Z(s),\alpha(s))ds\\
\ad   \le -\beta \int^{t_\varepsilon}_0
E_{z,\alpha}V(Z(s),\alpha(s))ds.\eea Given  $s
\le\tau_\varepsilon \wedge t$, we have $s\wedge \tau_\varepsilon=s$.
As a result, \bea\ad E_{z,\alpha}V(Z(t\wedge
\tau_\varepsilon),\alpha(t\wedge \tau_\varepsilon))-V(z,\alpha)\le
-\beta \int^t_0 E_{z,\alpha}V(Z(s\wedge
\tau_\varepsilon),\alpha(s\wedge\tau_\varepsilon))ds.\eea

By applying Gronwall's inequality, we obtain
 \bea\ad
E_{z,\alpha}V(Z(\tau_\varepsilon \wedge t), \alpha (\tau_\varepsilon
\wedge t)) \le V(z,\alpha)e^{-\beta t}. \eea Hence,\bea\ad
V_\varepsilon P(\tau_\varepsilon \le t)\le
E_{z,\alpha}[V(Z(\tau_\varepsilon), \alpha (\tau_\varepsilon)
)I_{\{\tau_\varepsilon \le t\}}] \le V(z,\alpha)e^{-\beta t},\eea in
which $V_\varepsilon=\inf\{V(z,i), z\in \rr^r \backslash
B_\varepsilon, i\in \mathcal {M}\}$ and $B_\varepsilon=\{z\in \K,
|z|< \varepsilon \}$, so $V_\varepsilon
>0$. Note that $\tau_\varepsilon \le t$ if and only if $\sup\limits_{0
\le u \le t}|Z(u)| \ge \varepsilon$. Therefore, it follows that
$P\{\sup\limits_{0 \le u \le t}|Z(u)| \ge \varepsilon \} \le \frac{
V(z,\alpha)e^{-\beta t}}{V_\varepsilon},$ so $P(|Z(t)|\ge
\varepsilon )\to 0$  as  $t\to \infty.$ That is,
$P(|X^{x,\alpha}(t)-X^{y,\alpha}(t)|\le \varepsilon)\to 1$ as
$t\rightarrow \infty.$ Thus, the proof is concluded. \qed

\subsection{State-Dependent Case}
Now, let us consider the case when the generator of the discrete
component $\alpha(t)$ is $x$ dependent. In this case, the switching
part is no  longer  a Markov chain. Because of the interplays
between $\alpha(t)$ and $X(t)$, we need more complex notations. We
use the same notations and  technique as that of \cite{BBG}.
Switching diffusions were treated in \cite{BBG}, whereas we deal
with switching jump diffusions. Define \beq{g}\barray\ad
\wdt{X}(t)=\left[X'(t)I_{\{\alpha(t)=1\}},
X'(t)I_{\{\alpha(t)=2\}},\cdots,X'(t)I_{\{\alpha(t)=m\}}\right]',\\
\ad
 S = \bigcup\limits_{i \in \cal{M}} 0_{r(i-1)} \times
\rr^r \times 0_{r(m-i)},\earray\eeq Here and in the sequel
$0_{k_1\times k_2}$ is a $\rr^{k_1 \times k_2}$ zero matrix , $0_k$
 denotes the $k$-dimensional zero
column vector. It is seen that
 $S\subseteq \mathbb{R}^{mr}$ and $\wdt{X}(t)$ is an $S$-valued
process. For $i\in\mathcal{M}, x\in \mathbb{R}^r$, define
 \beq{6.4} \barray\ad
 \tilde{x}^i= 0_{r(i-1)} \times x \times 0_{r(m-i)}\in S.\\
\ad \Xi =\bigcup\limits_{i,j \in \mathcal{M} \atop i<j}
0_{r(i-1)}\times \rr^r\times 0_{r(j-i-1)}\times \rr^r\times
0_{r(m-j)}. \earray\eeq
 Then $ \Xi\subseteq \mathbb{R}^{mr}$ and
$\wdt{X}^{x_0,i_0}(t)-\wdt{X}^{y_0,j_0}(t)$ is a $\Xi \cup S$-valued
process. For $x,y \in \mathbb{R}^r,i,j\in \mathcal{M}$, \bea\disp
\tilde x^i -\tilde y^j  = \left\{ \barray
 \ad [0'_{r(i - 1)},{x' - y'},0'_{r(m - i)}]' \in S  \text{   for  }\; i=j,\\
 \ad [0'_{r(i - 1)},{x'},0'_{r(j - i - 1)},- y', 0'_{r(m - j)}]'\in \Xi \text{  for  }\;i < j, \\
 \ad [0'_{r(j - 1)},{ - y'},
 0'_{r(i - j - 1)}, {x'}, 0'_{r(m - i)}]'
 \in\Xi \text{  for  }\,\;i > j.
 \earray \right.
\eea
Similar to the conditions we mentioned
in the previous part, under the condition {\rm(P1)} and (P2'), we
can obtain stability in distribution
similar to the approach in \cite{BBG}.
Now let us give condition (P2').

\begin{defn}\label{Pp}
{\rm The switching jump diffusion given by \eqref{2.1} and
\eqref{2.11} is said to satisfy condition (P2') if, for any
$\varepsilon
>0$ and any compact
subset $\K$ of $\rr^r$, there exists a $T=T(\varepsilon , \K)>0$ such
that \bea E|\wdt{X}^{x_0,i_0}(t)-\wdt{X}^{y_0,j_0}(t)|<\varepsilon
\text{ for all } t\ge T, \eea whenever $(x_0,i_0,y_0,j_0)\in \K
\times\mathcal{M} \times \K \times \mathcal{M}$. }\end{defn}

We can obtain (P2) from (P2'). To continue, we
focus on obtaining sufficient conditions for conditions (P1) and
(P2'). From \cite[Theorem 3.8]{BBG} we can see these two properties
imply asymptotic stability in distribution. So it is necessary to
establish sufficient criteria for the two properties. To proceed, we
need to introduce the following notations.

The generator $\mathcal {\tilde{G}}$ associated with the process
$\tilde x^i-\tilde y^j$ is defined as follows: For each $i,j \in
\mathcal {M}$, and for any twice continuously differentiable
function $f$, \bea \ad \mathcal {\tilde {G}} f(\tilde x^i-\tilde
y^j)=\mathcal {\tilde L}f(\tilde x^i-\tilde y^j)+\lambda
\int_\Gamma[f (\wdt{x}^i+\wdt {g}(x,i,\gamma)-\wdt{y}^j-\wdt{g}
(y,j,\gamma))-f(\tilde x^i-\tilde y^j)]\pi (d\gamma),\eea
 where $\mathcal{\tilde{L}}$ is the
 operator for a switching diffusion
process given by \beq{op-def} \barray
 {\cal \tilde L}f(\tilde x^i-\tilde y^j)\ad=\frac{1}{2}\tr(\tilde a(\tilde x^i,\tilde y^j)Hf(\tilde x^i-\tilde y^j))
 + (\tilde {b}(x,i)-\tilde {b}(y,j))'\nabla f(\tilde x^i-\tilde y^j)\\
 \ad + \sum^{m}_{k=1}q_{ik}(x)f(\tilde x^k-\tilde y^j)+
\sum^{m}_{k=1}q_{jk}(x)f(\tilde x^i-\tilde y^k)+\sum^{m}_{k=1 \atop
k\neq
i}\sum^{m}_{l=1 \atop l\neq j}\wdt
m(\Delta_{ik}(x)\cap\Delta_{jl}(y))\\
\aad \quad \times [f(\tilde x^k-\tilde y^l)-f(\tilde x^i-\tilde
y^l)-f(\tilde x^k-\tilde y^j)+f(\tilde x^i-\tilde y^j)],
 \earray\eeq
in which
 \bea\ad \tilde{b}(x,i)=\left[0'_{r(i-1)}, {b'(x,i)}, 0'_{r(m-i)}\right]',\\
 \ad \tilde{\sigma}(x,i)=\left[0'_{r(i-1)\times d},
 \sigma'(x,i), 0'_{r(m-i)\times d}\right]',\\
\ad \tilde{g}(x,i,\gamma)=\left[0'_{r(i-1)}, g'(x,i,\gamma), 0'_{r(m-i)}\right]',\\
\ad \tilde{a}(\tilde {x}^{i}, \tilde{y}^{j})=(\tilde{\sigma}(x,i)
-\tilde{\sigma}(y,j))
\times(\tilde{\sigma}(x,i)-\tilde{\sigma}(y,j))',\eea where
$0_{l_1\times l_2}$ is an $l_1\times l_2$ matrix with all entries
being 0, $b(x,i)$ and $g(x,i,\gamma) \in \rr^r$, and $\sigma(x,i)\in
\rr^{r\times d}$.   Recall that $\Delta_{ik}(x)$ are the intervals
having length $q_{ik}(x)$; $\wdt m$ is the Lebesgue measure on $\rr$
such that $dt \times \wdt m(dz)$ is the density of Poisson measure
with which we can represent the discrete component $\alpha(t)$ by a
stochastic integral as mentioned in Section \ref{sec:form}.

\begin{thm}
\label{Theorem 6.2} Assume the conditions of \thmref{Theorem 6.1}
hold and assume that for each $i,j\in\M$, there exists a Lyapunov
function $V(z)=z'z\in C^2(\rr^{mr}:\mathbb{R}_+)$ satisfying the
following condition:
 There exists a positive real
number $\wdt{\beta}$ such that \beq{cd-1} \mathcal {\tilde G}
V(\tilde x^i-\tilde y^j) \le -\wdt{\beta} V(\tilde x^i-\tilde
y^j),\eeq then \rm{(P1)} and \rm{(P2')} hold.
\end{thm}

\para{Proof.}  We need only verify (P2').
Let $\K$ be any compact subset of $\rr^r$, and fix any $x_0,y_0 \in
\K$, $i_0,j_0 \in \mathcal{M}.$ Define
\bea \ad \zeta_N =\inf\{t\ge
0,|\wdt{X}^{x_0,i_0}(t)-\wdt {X}^{y_0,j_0}(t)|>N\},\\
\ad
 \wdt\zeta_R=\inf\{t\ge 0,|\wdt {X}^{x_0,i_0}(t)|^2+|\wdt {X}^{y_0,j_0}(t)|^2>R\}.\eea
 Let $\zeta=\zeta_N \wedge \wdt \zeta_R$.

By virtue of the generalized It\^{o} formula, we have
  \bea
\disp E|\wdt{X}^{x_0,i_0}(t \wedge \zeta)-\wdt {X}^{y_0,j_0}(t\wedge
\zeta)|^2 =|\tilde x_0^{i_0}-\tilde y_0^{j_0}|^2+\int^{t\wedge
\zeta}_0 E\mathcal{\tilde G}|\wdt {X}^{x_0,i_0}(u)-\wdt
{X}^{y_0,j_0}(u)|^2 du.\eea Given the fact that for $u \le t\wedge
\zeta$, we have $u \wedge\zeta=u$. As a result, \bea \disp E|\wdt
{X}^{x_0,i_0}(t \wedge \zeta)-\wdt {X}^{y_0,j_0}(t\wedge \zeta)|^2
=|\tilde x_0^{i_0}-\tilde y_0^{j_0}|^2+\int^t_0 E\mathcal{\tilde
G}|\wdt {X}^{x_0,i_0}(u\wedge \zeta)-\wdt {X}^{y_0,j_0}(u\wedge
\zeta)|^2 du.\eea
 Then
 \bea \disp
 \frac{dE}{dt}|\wdt{X}^{x_0,i_0}(t \wedge
 \zeta)-\wdt
{X}^{y_0,j_0}(t\wedge \zeta)|^2\ad=E\mathcal{\tilde G}|\wdt
{X}^{x_0,i_0}(t\wedge \zeta)-\wdt {X}^{y_0,j_0}(t\wedge
\zeta)|^2\\
\ad \le -\wdt{\beta} E|\wdt {X}^{x_0,i_0}(t \wedge \zeta)-\wdt
{X}^{y_0,j_0}(t\wedge \zeta)|^2.\eea
Solving the differential inequality above leads to
 \bea E|\tilde{X}^ {x_0,i_0}(t \wedge
\zeta)-\tilde{X}^{y_0,j_0}(t \wedge \zeta)|^2\le e^{-\wdt{\beta}
t}|\tilde x_0^{i_0}-\tilde y_0^{j_0}|^2.\eea Let $N\to \infty, R\to
\infty$, we obtain  \bea E|\tilde{X}^
{x_0,i_0}(t)-\tilde{X}^{y_0,j_0}(t)|^2\le e^{-\wdt{\beta} t}|\tilde
x_0^{i_0}-\tilde y_0^{j_0}|^2.\eea Condition (P2') is thus verified.
\qed

\section{Examples}\label{sec:exm}
This section provides two simple examples. The purposes of these
examples are to demonstrate  our results.

\begin{exm}\label{6-1}{\rm
Consider a Markov switching jump diffusion
$\left(X(t),\alpha(t)\right)$ given by \eqref{2.1} and \eqref{2.11}
with $X(t)\in \rr^1$, $\lambda=1/ 8$,  $g(x,\alpha,\gamma)=x,$
$w(t)$ is a one-dimensional standard Brownian motion, $\alpha(t)\in
\mathcal{M}=\{1,2\}$ with \bea \ad Q = \left(
\begin{array}{ll}
   { - 1} & 1  \\
   3 & { - 3}  \\
\end{array} \right),
\  b(x,1)=\frac{x\sin x}{8},
\ b(x,2)=\frac{x\cos x}{2} ,
\
\sigma(x,1)=\frac{3x}{2},\ \sigma(x,2)=\frac{x}{2}.\eea Then by
choosing Lyapunov function $V=x^2e^{-4t}$, we can verify that
the equilibrium  $x=0$ is $p$-exponential stable with $p=2$.
It
is almost surely exponentially stable according to
\thmref{Theorem 4.2} as also confirmed in
\cite[Example 6.2]{Yinxi}. The following figure plots the
trajectories of the switched system.
\begin{figure}[h!]
\begin{center}
\includegraphics[width=8cm]{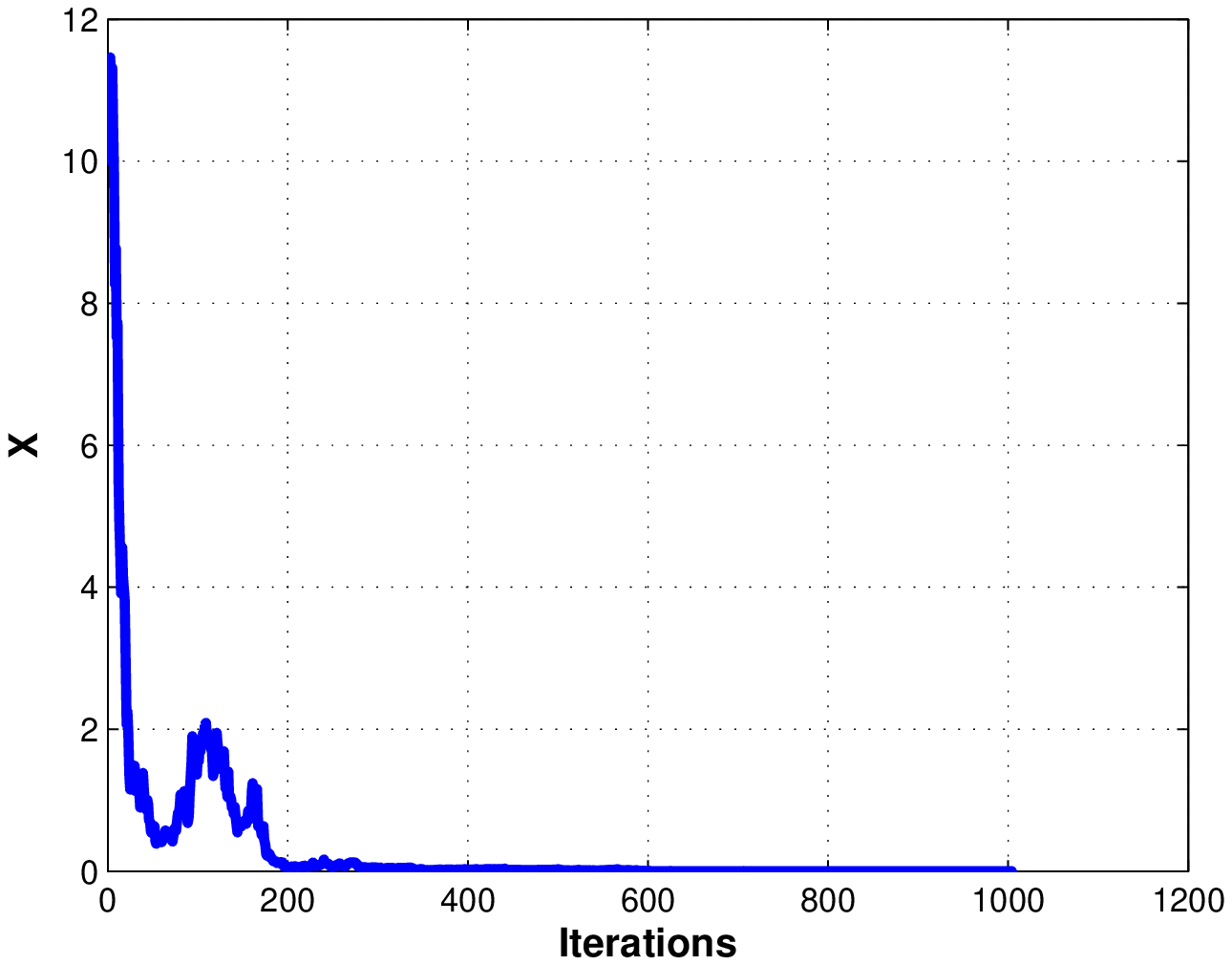}\\
\caption{Sample path of
\exmref{6-1} } \label{fig-1}
\end{center}
\end{figure}

}
\end{exm}

\begin{exm}\label{6-2}{\rm
Consider a switching jump diffusion $(X(t),\alpha(t))$ given  by
\eqref{2.1} and \eqref{2.11} with $X(t)\in \rr^1$, $\lambda=1$,
$g(x,\alpha,\gamma)=x,$ $w(t)$ is one-dimensional standard Brownian
motion, $\alpha(t)\in \mathcal{M}=\{1,2,3\}$ with
\[ Q(x) = \left( {\begin{array}{*{20}c}
   { - 3 - |\cos x| + \sin ^2 x \cos x} & {1 + |\cos x|} &
   {2 - \sin ^2 x \cos x}  \\
   1 & { - 1 - \frac{{x^2 }}{{1 + x^2 }}} & {\frac{{x^2 }}{{1 + x^2 }}}  \\
   {2 - \sin x\cos x} & {1 - \frac{{|x|}}{{1 + |x|}}\cos x} & { - 3 +
   \sin x \cos x + \frac{{|x|}}{{1 + |x|}}\cos x}  \\
\end{array}} \right),
\]
\bea \ad
b(x,1)=x+\sin x,\ b(x,2)=2x+x \sin x \cos x,
\ b(x,3)=3x+\sin^2 x,\\
\ad \sigma(x,1)=x+x\sin x, \  \sigma(x,2)=3x+x\cos x\sin x,\
\sigma(x,3)=x+\frac{x}{1+x}\sin x,\\
\ad g(x,i,\gamma)=x \text{ for }i=1,2,3. \ \hbox{ For the }  Q(x) \
\hbox{ given above}, \ \widehat{Q}= \left( {\begin{array}{*{20}c}
   { - 4} & 2 & 2  \\
   1 & { - 1} & 0  \\
   2 & 1 & { - 3}  \\
\end{array}} \right),\\
\ad
 b(1)=2, \ b(2)=2, \
b(3)=3,\ \sigma(1)=1,\ \sigma(2)=3, \ \sigma(3)=1.\eea Corresponding
to $\widehat{Q}$, the stationary distribution of the associated
Markov chain is given by
$\xi=(\frac{3}{13},\frac{8}{13},\frac{2}{13})$. There are three
associating jump diffusions that interact and switch back and forth.
They are given by \bea X(t)= x \ad + \int_0^t (X(s)+\sin X(s))ds
+\int_0^t( X(s)+X(s)\sin X(s))dw(s)\\
\ad   + \int_0^t {\int_\Gamma
{X(s^ - )} N(ds,d\gamma )}, \\
X(t)= x\ad +\int_0^t(2X(s)+X(s)\sin X(s)\cos X(s))ds\\
\ad + \int_0^t {(3X(s)+X(s)\cos X(s)\sin X(s))dw(s)}
+ \int_0^t {\int_\Gamma {X(s^
- )} N(ds,d\gamma )}, \\
 X(t)= x\ad +\int_0^t (3X(s)+\sin^2 X(s))ds
+\int_0^t( X(s)+\frac{X(s)}{1+X(s)}\sin X(s))dw(s)\\
\ad  + \int_0^t {\int_\Gamma {X(s^ - )} N(ds,d\gamma )}. \eea
From \cite[Example 6.1]{Yinxi}, the first and the third jump
diffusion are unstable in probability and therefore, not
asymptotically stable in the large. In view of \cite{Wee} and
\thmref{thm:st-l}, the second jump diffusion is asymptotically stable
in the large. Using \corref{cor-34},
 we obtain that the switching jump
 diffusion is asymptotically stable in
the large. The following plot provides a  sample
path of the regime-switching jump diffusion.
\begin{figure}[h!]
\begin{center}
\epsfig{figure=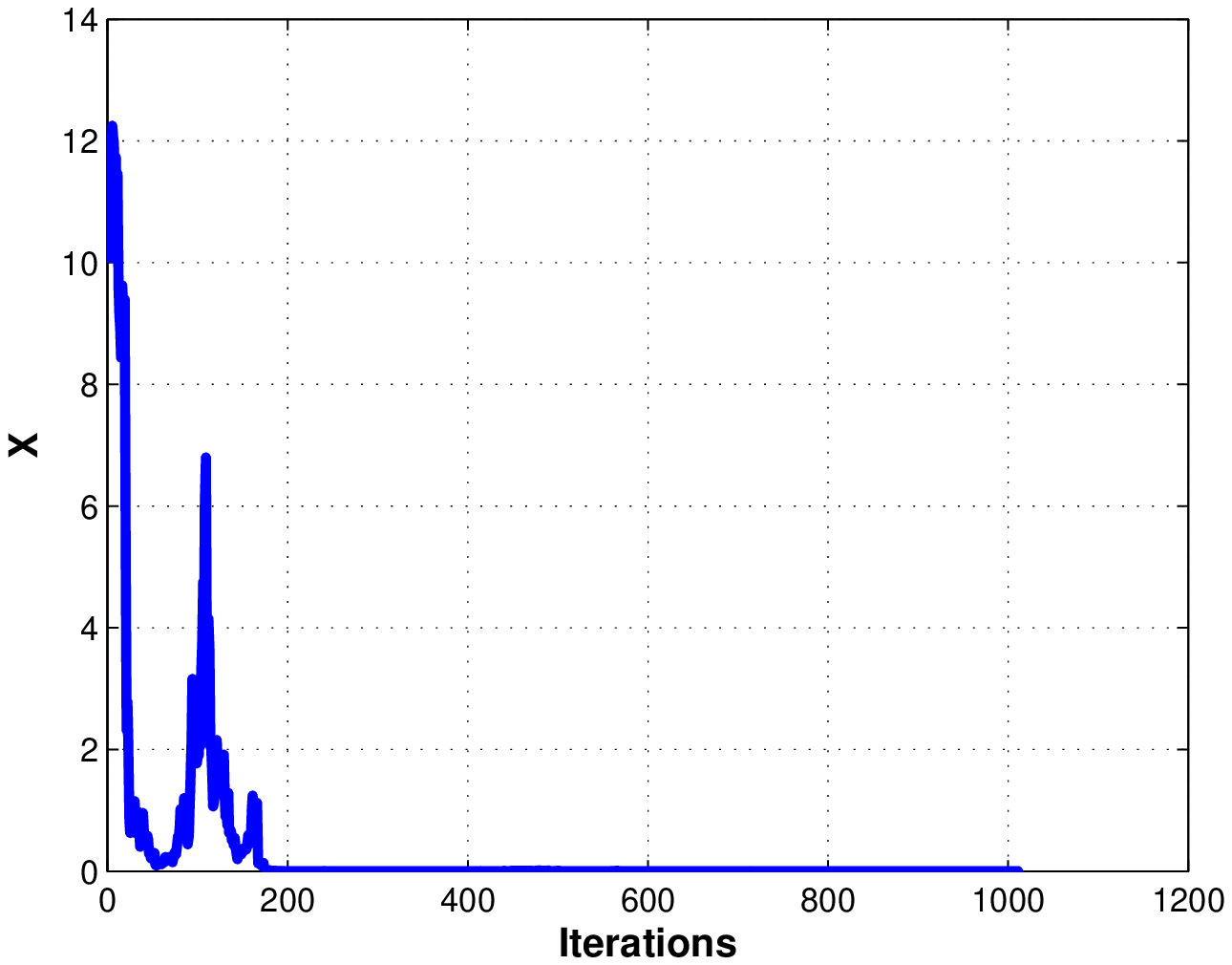,width=0.5\linewidth }
\caption{Sample path of
\exmref{6-2}}\label{fig-2}
\end{center}
\end{figure}

}
\end{exm}

\section{Further Remarks}\label{sec:fur}
 This work focused on stability of regime
switching jump diffusions.
Under simple conditions, we derived
  sufficient conditions
for asymptotic stability in the large and asymptotic stability
in distribution. We also provided necessary and sufficient
conditions for exponential stability. The connection between
exponential stability and almost surely exponential stability was
studied. Smooth dependence on the initial data was demonstrated as
well. Future research efforts can be directed to the study of
positive recurrence  and egrodicity of regime-switching jump
diffusions, which was coined as weak stability in \cite{Wonham}
for diffusion processes. Treating stability with non-Lipschitz
coefficients is of great importance.  Numerical methods will be
welcomed since the nonlinear systems rarely have closed-form
solutions. All of these deserve more thoughts and further
considerations.

\end{document}